\documentclass[10pt]{article}
\usepackage{amsmath,amssymb,wasysym}
\usepackage{graphicx}
\usepackage{mathrsfs}
\usepackage{pst-all}
\usepackage{nonfloat}
\usepackage[left=3cm,right=3cm, top=3.3cm, bottom=3cm]{geometry}
\newtheorem{theo}{Theorem}[section]
\newtheorem{lemm}[theo]{Lemma}

\newtheorem{prop}[theo]{Proposition}

\newtheorem{remark}[theo]{Remark}
\def\proof {{\noindent \bf{Proof:\hspace{4pt}}}}
\def\endproof{\hfill$\square$\vspace{6pt}}
\numberwithin{equation}{section}
\title{
{\bf\Large  Time-decay and Strichartz estimates for the Benjamin-Bona-Mahony equation and existence of solutions on modulation spaces}}
\author{
{\bf\large Carlos Banquet}\footnote{Email: cbanquet@correo.unicordoba.edu.co}\\
{\it\small  Departamento de Matem\'{a}ticas y Estad\'{\i}stica}\\ 
{\it\small Universidad de C\'{o}rdoba}\\
 {\it\small A.A. 354, Monter\'{\i}a, Colombia}\\
 \vspace{0mm}\\
 {\bf\large \'Elder  J. Villamizar-Roa}\footnote{Email: jvillami@uis.edu.co}\hspace{2mm}
{\bf\large}\\
{\it\small Escuela de Matem\'aticas}\\
 {\it\small Universidad Industrial de Santander}\\
{\it\small  A.A. 678, Bucaramanga, Colombia}\vspace{3mm}\\}
\date{}
\begin{document}
\maketitle
\begin{abstract}
In this paper we derive time-decay and Strichartz estimates for the generalized Benjamin-Bona-Mahony equation on the framework of modulation spaces $M^s_{p,q}.$ We use this results to analyze the existence of local and global solutions of the corresponding Cauchy problem with rough data in modulation spaces. The results improve known results in Sobolev spaces in some sense. 

{\bf Key words.} Benjamin-Bona-Mahony equation, modulation spaces, Strichartz estimates, well-posedness.\\

{\bf AMS subject classifications.} 35Q53; 35A01; 35Q35; 35C15
\end{abstract}

\section{Introduction}\label{intro}
 In this paper we study the Cauchy problem for the generalized Benjamin-Bona-Mahony equation, (gBBM)
\begin{equation}\label{BBMequa}
\left\{
\begin{array}{lc}
u_{t} +u_x-u_{xxt}+u^{\lambda}u_x=0, & x\in \mathbb{R},\ \ t\in \mathbb{R},
\\
u(x,0)=u_{0}(x), \ \ \ & x\in \mathbb{R},
\end{array}
\right.  
\end{equation}
where $u: \mathbb{R} \times \mathbb{R} \rightarrow \mathbb{R}$ is a real-valued function, $u_0:\mathbb{R} \rightarrow \mathbb{R}$ is the initial data and $\lambda\geq 1$ is an integer. The case $\lambda=1$ corresponds to the Benjamin-Bona-Mahony equation (BBM), which has been derived as a model to describe the gravity water waves in the long-wave
 regime, see Benjamin, Bona and Mahony \cite{BenBonMah}, and Peregrine \cite{Peregrine1, Peregrine2}. Also, BBM equation is well suited for modeling wave propagation on star graphs, which gives some interesting applications as shown in Bona and Cascaval \cite{BonCas}. For $\lambda=2$ the equation (\ref{BBMequa}) is known as the modified BBM equation, which describes wave propagation in one dimensional nonlinear lattice (cf. Wadati \cite{Wadati1, Wadati2}); thus, the generalization considered in (\ref{BBMequa}) is not only of mathematical interest.\\

The BBM equation is a good substitute for the famous Korteweg-de Vries equation (KdV) 
\[u_t+u_x-u_{xxx}+u_xu=0, \ \ \ x,t \in\mathbb{R},\]
in the case of shallow waters in a channel (see Whitham \cite{Whitham}, Bona, Pritchard and Scott \cite{BonPriSco}). Furthermore, the solutions of the KdV and the BBM  stay ``close'' to each other over relatively long time intervals, see \cite{BonPriSco} for more details.   As the KdV equation, the gBBM possesses solitary and periodic wave solutions, which are particular solutions very important for applications in physics. The existence, orbital, asymptotic and spectral stability or instability of the solitary or periodic traveling waves, associated to the gBBM equation, have been studied by several researchers, see for instance \cite{AngBanSci1, AngBanSci2, eldika, haragus1, MillerWeinstein, SougaStrauss, weinstain1, zeng}.\\

The mathematical analysis of well-posedness and ill-posedness of (\ref{BBMequa}) has been considered extensively in the literature (see \cite{Bona, Carvajal, Panthee, Roum, Wang2, Ming_Wang} and references therein). The well-posedness of (\ref{BBMequa}) with $\lambda=1,$ in Sobolev spaces $H^s(\mathbb{R}),$ $s\geq 0,$ was obtained by Bona and Tzvetkov in \cite{Bona}, and Carvajal and Panthee \cite{Carvajal, Panthee} in the periodic case on $H^s_{per}([-L,L])$ with $s\geq 0$. On the other hand, the initial value problem (\ref{BBMequa}) with $\lambda=1$ is ill-posed in $H^s(\mathbb{R})$ or $H^s_{per}([-L,L]),$ for $s<0$ (cf. \cite{Bona, Panthee}). Previous results have been obtained in finite energy spaces $H^s.$ Results of well-posedness for a generalization on the space dimension of the BBM equation, in Sobolev spaces $W^{k,p},$ have been obtained by Goldstein and Wichnoski  \cite{Goldstein}, and Avrin and Goldstein \cite{Avrin}. In \cite{Goldstein} the authors analyzed the local well-posedness in $W^{2,p}(\Omega)\cap W^{1,p}_0(\Omega)$ for $p>n,$ with $\Omega$ a bounded domain of $\mathbb{R}^n,$ and in \cite{Avrin} the authors studied the existence of local weak solutions in $W^{1,p}(\mathbb{R}^n),$ $p\geq 1.$ Recently, Wang \cite{Ming_Wang} established an interesting result of global well-posedness for the BBM equation in Bessel potential spaces $H_p^s(\mathbb{R}),$ with $s\geq\ \mbox{max}\{0,\frac{1}{p}-\frac{1}{2}\}$ and $1\leq p<\infty.$  The results of \cite{Ming_Wang} are sharp in the sense that equation (\ref{BBMequa}) with $\lambda=1$ is ill-posed in $H^{s}_p(\mathbb{R})$ for $s<\ \mbox{max}\{0,\frac{1}{p}-\frac{1}{2}\},$ $1\leq p<\infty.$ More recently, Bona and Dai in \cite{BonDai} obtained an ill-posedness result for the BBM equation on the periodic homogeneous Sobolev spaces $\dot{H}^r_{per},$ with $r<0$. So far, to the best of our knowledge, the larger initial data classes for BBM equation are those of \cite{Ming_Wang}.\\

Our interest in this paper is to analyze the well-posedness in some spaces of low regularity than $H^s$ and $H^s_p$ for large $s,$ namely, modulation spaces $M^s_{p,q}. $ Modulation spaces are decomposition spaces that emerge from a uniform covering of the underlying frequency space; they were introduced by Feichtinger in \cite{Feichtinger}, prompted by the idea of measuring the smoothness classes of functions or distributions.  Since their introduction, modulation spaces have become canonical for both time-frequency and phase-space analysis, see Chaichenets {\it et al.} \cite{Chai}. Wang and Hudzik \cite{BaoxHudz} gave an equivalent definition of modulation spaces by using the frequency-uniform-decomposition operators. In the same work, the existence of global solutions for nonlinear Schr\"odinger and Klein-Gordon equations in modulation spaces were analyzed. After them, several studies on nonlinear PDEs in the framework of modulation spaces have been addressed (cf. \cite{Chai,Huang,Iwabuchi,Manna,Ruz,Wang,Zhao} and references therein). In this context, the contribution of this paper is to analyze the existence of solutions for the gBBM equation with initial data in modulation spaces. To get this aim, first we establish a careful harmonic analysis in order to derive some time-decay estimates of the one parameter group given by the corresponding linear equation, as well as some Strichartz estimates and nonlinear estimates on modulation spaces which allow us to control the nonlinearity in the gBBM equation (cf. Section \ref{Section2} and \ref{Section3}). In particular, we prove some Strichartz estimates in $M^{s}_{p,q}$ for a general dispersive semigroup $U(t)=\mathscr{F}^{-1}e^{itP(\xi)}\mathscr{F},$ with $P:\mathbb{R}^n\to\mathbb{R}$ a real-valued function, complementing the ones established in Wang and Hudzik \cite{BaoxHudz}, which can be used to analyze the well-posedness of another dispersive models.\\ 

Before stating our main results, we recall some preliminar definitions and notations related to the modulation spaces $M_{p,q}^s(\mathbb{R}^n)$ (for more details see for instance Wang and Hudzik \cite{BaoxHudz} and Kato \cite{Kato}). Let $\mathscr{S}$ be the Schwartz space on $\mathbb{R}^n$ and $\mathscr{S}^\prime$ its dual space. Let $Q_0=\{\xi:\xi_i\in [1/2,1/2),\ i=1,...,n\}$ and $Q_k=k+Q_0,$ $k\in \mathbb{Z}^n.$ Thus, $\{Q_k\}_{k\in \mathbb{Z}^n}$ constitutes a decomposition of $\mathbb{R}^n,$ that is, $Q_i\cap Q_j=\emptyset$ and $\bigcup_{k\in \mathbb{Z}^n}Q_k=\mathbb{R}^n.$ Let $\rho:\mathbb{R}^n\rightarrow [0,1]$ be a smooth function satisfying $\rho(\xi)=1$ for $\vert \xi\vert\leq \frac{\sqrt{n}}{2}$ and $\rho(\xi)=0$ for $\vert \xi\vert\geq \sqrt{n}.$ Let $\rho_k(\xi)=\rho(\xi-k),$ $k\in\mathbb{Z}^n,$ a translation of $\rho.$ It holds that $\rho_k(\xi)=1$ in $Q_k,$ and thus, $\sum_{k\in\mathbb{Z}^n}\rho_k(\xi)\geq 1$ for all $\xi\in \mathbb{R}^n.$ Let
$$\sigma_k(\xi)=\rho_k(\xi)\left( \sum_{l\in\mathbb{Z}^n}\rho_l(\xi)\right)^{-1},\ k\in\mathbb{Z}^n.$$
Then, the sequence $\{\sigma_k(\xi)\}_{k\in\mathbb{Z}^n}$ verifies the following properties:
\begin{eqnarray*}
&\vert \sigma_k(\xi)\vert\geq C,\ \forall\xi\in Q_k,&\\
&\mbox{supp}(\sigma_k)\subset\{\xi: \vert \xi-k\vert_\infty\leq \sqrt{n}\},&\\
&\sum\limits_{k\in\mathbb{Z}^n}\sigma_k(\xi)=1,\ \forall\xi\in\mathbb{R}^n,&\\
&\vert D^\alpha\sigma_k(\xi)\vert\leq C_{m},\ \forall\xi\in \mathbb{R}^n,\ \vert\alpha\vert\leq m.&
\end{eqnarray*}
Modulation spaces $M^s_{p,q}=M^s_{p,q}(\mathbb{R}^n)$ are Banach spaces constituted by frequency uniform decomposition $\sigma_k,$ $k\in \mathbb{Z}^n.$ Explicitly, we consider the frequency-uniform decomposition operators $\square_k:=\mathscr{F}^{-1}\sigma_k\mathscr{F}=\mathscr{F}^{-1}[\sigma_k\cdot\mathscr{F}],$ $k\in\mathbb{Z}^n.$ Then, for $s\in\mathbb{R},$ $1\leq p,q\leq\infty,$ modulations spaces $M_{p,q}^s$ are defined as (cf. \cite{Feichtinger, Kato, BaoxHudz}): 
\begin{eqnarray*}
M_{p,q}^s:=\left\{ f\in \mathscr{S}^{\prime}(\mathbb{R}^n):\ \Vert f\Vert_{M^s_{p,q}}<\infty\right\},
\end{eqnarray*}
where
\begin{equation*}
\Vert f\Vert_{M^s_{p,q}}=\left\{
\begin{array}{lc}
 \left( \sum\limits_{k\in\mathbb{Z}^n}(1+\vert k\vert)^{sq}\Vert \square_kf\Vert_p^q\right)^{1/q},\ \ \mbox{for}\ 1\leq q<\infty,\\
 \\
\sup\limits_{k\in\mathbb{Z}^n}(1+\vert k\vert)^{s}\Vert \square_kf\Vert_p, \ \ \mbox{for}\ q=\infty.
\end{array}
\right.
\end{equation*}
For simplicity, we will write $M^0_{p,q}(\mathbb{R})=M_{p,q}(\mathbb{R}).$ Many of their properties, including embeddings
in other known function spaces, can be found in Wang and Hudzik  \cite{BaoxHudz} (see also Kato \cite{Kato}). In particular, the following properties hold:
\begin{itemize}
\item[i)]If $\Omega$ is a compact subset of $\mathbb{R}^n,$ then $\mathscr{S}^\Omega=\{f:f\in \mathscr{S}\ \mbox{and}\ \mbox{supp}\widehat{f}\subset \Omega\}$ is dense in $M^s_{p,q},$ $s\in\mathbb{R},$ $0<p,q<\infty.$
\item[ii)] $M^{s_1}_{p_1,q_1}\subset M^{s_2}_{p_2,q_2},\ \mbox{if}\ s_1\geq s_2,\ 0<p_1\leq p_2,\ 0<q_1\leq q_2.$
\item [iii)] $M^{s_1}_{p,q_1}\subset M^{s_2}_{p,q_2},\ \mbox{if}\ q_1>q_2,\ s_1>s_2,\ s_1-s_2>n/q_2-n/q_1.$
\item [iv)]  $M_{p,1}\subset L^\infty\cap L^p,$  for $1<p\leq \infty.$ 
\item [v)]  $B^{s+n/q}_{p,q}\subset M^s_{p,q},$ for $0<p,q\leq \infty$ and $s\in\mathbb{R}.$ 
\item [vi)]  $B^{s_1}_{p,q}\subset M^{s_2}_{p,q}$ if and only if $s_1\geq s_2+n\nu_1(p,q).$
\item [vii)]  $H^{s_1}_{p}\subset M^{s_2}_{p,q},$ if  $s_1> s_2+n\nu_1(p,q),$ where 
\begin{equation*}
\nu_1(p,q)=\left\{
\begin{array}{lc}
0,\ \ \mbox{if}\ (\frac{1}{p},\frac{1}{q})\in \left\{(\frac{1}{p},\frac{1}{q})\in [0,\infty)^2\ \ :\ \frac{1}{q}\leq \frac{1}{p}\ \mbox{and}\ \frac{1}{q}\leq 1-\frac{1}{p}\right\},\\
\\
\frac{1}{p}+\frac{1}{q}-1,\ \ \mbox{if}\ (\frac{1}{p},\frac{1}{q})\in \left\{(\frac{1}{p},\frac{1}{q})\in [0,\infty)^2\ \ :\ \frac{1}{p}\geq\frac{1}{2}\ \mbox{and}\ \frac{1}{q}\geq 1-\frac{1}{p}\right\}, \\
\\
-\frac{1}{p}+\frac{1}{q},\ \ \mbox{if}\ (\frac{1}{p},\frac{1}{q})\in \left\{(\frac{1}{p},\frac{1}{q})\in [0,\infty)^2\ \ :\ \frac{1}{p}\leq\frac{1}{2}\ \mbox{and}\ \frac{1}{q}\geq \frac{1}{p}\right\}.
\end{array}
\right.
\end{equation*}
\end{itemize}
%
In order to establish the main results, we need to consider the integral formulation associated to the Cauchy problem (\ref{BBMequa}). Applying the operator $(1-\partial_{xx})^{-1}$ on both sides of  (\ref{BBMequa}) it holds that 
\begin{equation}
\left\{
\begin{array}{lc}
u_{t} =i\varphi(D)u-\frac{i}{2}\varphi(D)[u^{\lambda+1}], & x\in \mathbb{R},\ \ t\in \mathbb{R},
\\
u(x,0)=u_{0}(x), \ \ \ & x\in \mathbb{R},
\end{array}
\right.  \label{BBMequa2}
\end{equation}
where $\varphi(D)$ is defined as the Fourier multiplier with symbol $\varphi(\xi)=\frac{\xi}{1+\xi^2}.$ Let $S(t)$ be the unitary group on $L^2$ generated by $-i\varphi(D),$ namely, $S(t)u_0=e^{-it\varphi(D)}u_0$ with $\widehat{\varphi(D)u_0}(\xi)=\varphi(\xi)\widehat{u_0}(\xi).$ Then, by the Duhamel principle (\ref{BBMequa2}) is equivalent to the
following  integral equation
\begin{equation}
u(x,t)=S(t)u_{0}(x)-\frac{i}{\lambda+1}\int_{0}^{t}S(t-\tau )\varphi(D)[u^{\lambda+1}(x,\tau)] d\tau.
\label{IntEqu}
\end{equation}
Now we are in position to establish the main results of this paper. From now on, we consider the function $\beta_\sigma:(-\infty,-1)\longrightarrow [-\frac{1}{3},0)$ defined by:
\begin{equation}\label{Beta}
\beta_{\sigma}:=\beta(\sigma)=\left\{
\begin{array}{cc}
\frac{\sigma+1}{1-2\sigma}, & \mbox{if} \ \ -4\leq \sigma<-1,\vspace{0.2cm}\\
\frac{-3}{1-2\sigma}, & \mbox{if} \ \ -\infty <\sigma\leq -4.
\end{array}
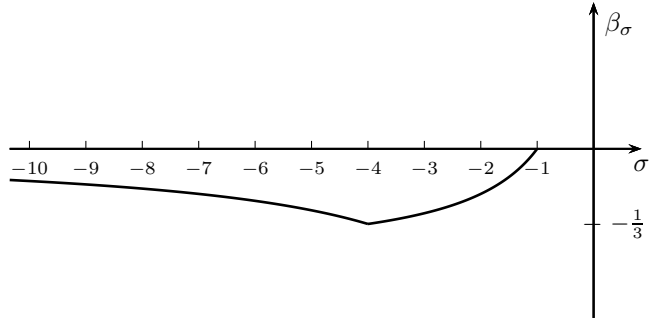
\right. 
\end{equation}

\begin{minipage}{\linewidth}
\begin{center}
\psset{xunit=0.75cm,yunit=3.0cm}
\begin{pspicture}(-10.35,-0.65)(0.85,0.65)
\psaxes[labelFontSize=\scriptstyle,ticksize=3pt,arrowsize=0.125,linewidth=0.7pt]{->}(0,0)(-10.35,-0.75)(0.85,0.65)
\psplot[plotstyle=curve,linewidth=1pt]
       {-10.35}{-4}{3 1 2 x mul neg add div neg} %
\psplot[plotstyle=curve,linewidth=1pt]
       {-4}{-1}{1 x add 1 2 x mul neg add div} %
\rput[l](0.2,0.55){$\beta_{\sigma}$}
\rput[b](0.85,-0.1){$\sigma$}
\rput[l](0.3,-0.333){\small$-\frac13$}
\rput[l](-0.2,-0.333){$-$}
\end{pspicture}\\
\end{center}
\end{minipage}
\vspace{0.4cm}
\\
Function $\beta_{\sigma}$ controls the time-decay of the $L^{p'}-H^\sigma_p, \sigma<-1,$ estimate of the group $S(t),$ which plays a key rol in the analysis of the existence of global solutions for (\ref{IntEqu}).
\begin{theo}\label{GlobalLamb0}
Consider $\lambda\geq 1$ an integer, $p=\lambda+2,$  $s\geq 0,$ and $0<\theta\leq -\frac{1}{\sigma},$ with $\sigma<-1$ such that $0<\frac{-\lambda\theta\beta_{\sigma}}{\lambda+2}<1$ and define $r=\frac{\lambda(\lambda+2)}{\lambda+2+\lambda\theta\beta_{\sigma}}.$ Then, there exists $\epsilon>0$ such that if $\Vert S(t)u_0\Vert_{L^r(\mathbb{R};M^s_{p,1})}<\epsilon,$ equation (\ref{IntEqu}) has a unique global solution $ u\in L^{r}(\mathbb{R}; M^s_{p,1}).$
\end{theo}
The restriction $q=1$ in Theorem \ref{GlobalLamb0} comes from the use of the product estimate in modulation spaces (cf. Lemma \ref{prod} below). Considering a different product estimate given in Iwabuchi \cite{Iwabuchi} we can include a global existence result with values in $M^s_{p,q},$ for $1\leq q<2$ by penalizing the regularity coefficient $s,$ imposing that $1-\frac 1q\leq s<\frac 1q.$ This is the content of next theorem.
\begin{theo}\label{GlobalLamb1}
Consider $\lambda\geq 1$ an integer, $p=\lambda+2,$  $1\leq q <\infty $  and $1-\frac 1q\leq s<\frac 1q.$ Also let $0<\theta\leq -\frac{1}{\sigma},$ with $\sigma<-1$ such that $0<\frac{-\lambda\theta\beta_{\sigma}}{\lambda+2}<1$ and define $r=\frac{\lambda(\lambda+2)}{\lambda+2+\lambda\theta\beta_{\sigma}}.$ Then, there exists $\epsilon>0$ such that if $\Vert S(t)u_0\Vert_{L^r(\mathbb{R};M^s_{p,q})}<\epsilon,$ equation (\ref{IntEqu}) has a unique global solution $ u\in L^{r}(\mathbb{R}; M^s_{p,q}).$
\end{theo}
In Theorems \ref{GlobalLamb0} and \ref{GlobalLamb1}, the condition $0<\frac{-\lambda\theta\beta_{\sigma}}{\lambda+2}<1$ and the value $r=\frac{\lambda(\lambda+2)}{\lambda+2+\lambda\theta\beta_{\sigma}}$ come from the application of the Hardy-littlewood-Sobolev's inequality, which allow us to consider the range $\lambda\geq 1.$ However, we are forced to assume a smallness condition on the norm $\Vert S(t)u_0\Vert_{L^r(\mathbb{R};M^s_{p,q})},$ in place of a weaker smallness condition on the initial data $u_0$ directly. By using Strichartz estimates we can control the nonlinearity and establish the existence of global solution by assuming $\Vert u_0\Vert_{M^{s-\frac{\sigma\theta}{2}}_{2,q}}$  small enough; in this case, we need to impose the condition $\lambda\geq 6.$ This is the content of next theorem.
 \begin{theo}\label{TheoX0}
Consider $\lambda\geq 6$ an integer, $0<\theta\leq -\frac{1}{\sigma},$ with $\sigma=\frac{\lambda+2}{4-\lambda}$ or $\sigma=\frac{2-3\lambda}{4}.$ Take $p=\lambda+2,$  $q\in [\gamma',\gamma]$ with $\gamma=-\frac{2(\lambda+2)}{\lambda\theta\beta_{\sigma}}$ and $1-\frac 1q\leq s<\frac 1q.$ Then, there exists $\epsilon>0$ such that if $\Vert u_0\Vert_{M^{s-\frac{\sigma\theta}{2}}_{2,q}}<\epsilon,$ equation (\ref{IntEqu}) has a unique global solution $u$ in   
$C(\mathbb{R}; M^s_{2,q}) \cap L^{\gamma}(\mathbb{R}; M^s_{p,q}).$
\end{theo}
By using only the time-decay estimates for the group $S(t)$ on Modulation spaces, we are able to obtain the existence of global solutions in the time-weighted based space $X_{p,q}^{\rho,s}$ based on the modulation spaces:
$$X_{p,q}^{\rho,s}:=\left\{u\in L^\infty_{loc}(\mathbb{R};M^s_{p,q}):\ \sup_{-\infty<t<\infty}(1+\vert t\vert)^{\rho}\Vert u(t)\Vert_{M^s_{p,q}}<\infty\},\ \rho=-2\theta\left(\frac{1}{2}-\frac{1}{p}\right)\beta_\sigma>0\right\}.$$ 
 \begin{theo}\label{TheoX2a}
Consider $\lambda\geq 3$ an integer, $0<\theta\leq -\frac{1}{\sigma},$ with $\sigma<-1,$  $p=\lambda+2,$ $s\geq 0.$  Then, there exists $\epsilon>0$ such that if $\Vert u_0\Vert_{M^{s-{\sigma\theta}}_{p',1}}<\epsilon,$ equation (\ref{IntEqu}) has a unique global solution $u\in X_{p,1}^{\rho,s}.$ 
\end{theo}
 \begin{theo}\label{TheoX2}
Consider $\lambda\geq 3$ an integer, $0<\theta\leq -\frac{1}{\sigma},$ with $\sigma<-1,$  $p=\lambda+2,$ $1\leq q<\infty,$ and $1-\frac 1q\leq s<\frac 1q.$ Then, there exists $\epsilon>0$ such that if $\Vert u_0\Vert_{M^{s-{\sigma\theta}}_{p',q}}<\epsilon,$ equation (\ref{IntEqu}) has a unique global solution $u\in X_{p,q}^{\rho,s}.$
\end{theo}
In Theorems \ref{TheoX2a} and \ref{TheoX2}, the condition $\lambda\geq 3$ comes from the integrability of the function $f(\tau)=(1+\vert t-\tau\vert)^{-\rho}(1+\vert \tau\vert)^{-\rho(\lambda+1)},$ on the real line. Next theorems provide local existence results in $C([-T,T];M^{s}_{p,q})$ for the general case $\lambda\geq 1.$

 \begin{theo}\label{TheoX10}
Consider $\lambda\geq 1$ an integer, $2\leq p<\infty,$ $s\geq 0,$  and assume $u_0\in M^s_{p',q}.$ Then there exists $T>0$ and a unique solution $u\in C([0,T];M^{s}_{p,1})$ solution of equation (\ref{IntEqu}).
 \end{theo}
  \begin{theo}\label{TheoX10b}
Consider $\lambda\geq 1$ an integer, $2\leq p<\infty,$ $1\leq q <\infty,$ $1-\frac 1q\leq s<\frac 1q,$  and assume $u_0\in M^s_{p',q}.$ Then there exists $T>0$ and a unique solution $u\in C([0,T];M^{s}_{p,q})$ solution of equation (\ref{IntEqu}).
 \end{theo}
 \begin{remark}
\begin{itemize}
\item [i)]  Theorems \ref{GlobalLamb0}, \ref{GlobalLamb1}, \ref{TheoX0}, \ref{TheoX2a} and \ref{TheoX2}, continue true if we replace the time interval $\mathbb{R}$ by the compact interval $[-T,T]$ throughout their statements.  Notice that for $s_1>s+\nu_1(p,q),$ $p=\lambda+2,$ $H^{s_1}_p\subset M^s_{p,q},$ and therefore for $1<r\leq \infty$ it holds
\begin{eqnarray*}
\Vert S(t)u_0\Vert_{L^r([-T,T];M^s_{p,q})}\leq C \Vert S(t)u_0\Vert_{L^r([-T,T];H^s_{p})}\leq CT\Vert u_0\Vert_{H^s_p}.
\end{eqnarray*}
In particular, observing the statement of Theorem \ref{GlobalLamb1} for instance,  if $u_0\in H^s_p$ and $0<T\leq \epsilon C^{-1}\Vert u_0\Vert^{-1}_{H^s_p},$ where $\epsilon$ is as in Theorem  \ref{GlobalLamb1}, it follows that $\Vert S(t)u_0\Vert_{L^r([-T,T];M^s_{p,q})}\leq \epsilon.$ Thus, without assume any smallness condition on the initial data, the local in time version of Theorem  \ref{GlobalLamb1} gives a solution $u\in L^r([-T,T];M^s_{p,q}),$ for initial data in $H^s_p.$ Therefore, if $\tilde{u}\in C([-T_0,T_0];H^s_p)$ is the local solution obtained in Wang \cite{Ming_Wang}, for some $T_0>0,$ by uniqueness $u=\tilde{u}.$ 
\item [ii)] The initial data class in Theorems \ref{TheoX0}, \ref{TheoX2a} and \ref{TheoX2} is 
larger than the $H^{s_1}_p$ of Wang \cite{Ming_Wang} provided $s_1$ be large enough. This is consequence of the embedding $H^{s_1}_{p}\subset M^{s_2}_{p,q},$ if  $s_1> s_2+n\nu_1(p,q).$
\item [iii)] In Theorems \ref{GlobalLamb1}, \ref{TheoX0}, \ref{TheoX2} and \ref{TheoX10b}, the condition $1-\frac 1q\leq s<\frac 1q$ and the integer nature of $\lambda$ come from Lemma \ref{Ib2} below, to use the estimate $\|u^p\|_{M^{s}_{{q},\mu}}\leq C  \|u\|^p_{M^{s}_{{pq},\nu}}.$ In fact, it is an open problem to see if $\|u^p\|_{M^{s}_{{q},\mu}}\leq C  \|u\|^p_{M^{s}_{{pq},\nu}}$ holds for any positive real constant $p.$
\end{itemize}
 \end{remark}
 This paper is organized as follows. In Section 2, we state time-decay  and the Strichartz estimates in $M^s_{p,q}$ for the group $S(t),$ as well as some nonlinear estimates to deal with the nonlinear term in gBBM equation. In Section 3, we prove some Strichartz estimates in $M^{s}_{p,q}$ for a general dispersive semigroup $U(t)=\mathscr{F}^{-1}e^{itP(\xi)}\mathscr{F},$ with $P:\mathbb{R}^n\to\mathbb{R}$  a real-valued function, which, applied to the particular case $S(t),$ allow us to obtain some existence results. Finally, we prove Theorems \ref{GlobalLamb0}-\ref{TheoX10b} in Section 4.
\section{Dispersive and nonlinear estimates}\label{Section2}
The first aim of this section is to derive some decay estimate of the group $S(t)$ on Modulation spaces $M^s_{p,q}.$ For that, a useful tool is the van der Corput's Lemma, whose proof can be found in Stein \cite{Stein}, see also Linares and Ponce \cite{LinaresPonce}.
\begin{lemm}[van der Corput]\label{VanCorput}
Let $h$ be either convex or concave twice differentiable function and $F$  be continuously differentiable function on $[a, b],$ with $-\infty\leq a<b\leq \infty.$  Then
\begin{eqnarray}
\left\vert \int_a^b F(r)e^{ih(r)}dr\right\vert \leq 4\left\{ \min_{r\in [a,b]}\left\vert h''(r)\right\vert\right\}^{-1/2}\left [|F(b)|+\int_a^b\vert F'(r)\vert dr\right],
\end{eqnarray}
for $h''\neq 0$ on $[a,b].$
\end{lemm}
Using Lemma \ref{VanCorput} and a meticulous analysis of the Fourier symbol of the $S(t)$ we can obtain the following $L^1-L^\infty$ estimate.
\begin{lemm}\label{EstEpsEne}
Let  $\sigma<0,$ $0<\epsilon<\frac{1}{8}$ and $N>2.$  Define $K(t)=S(t)J^{\sigma},$ where $J^\sigma=(I-\Delta)^{\sigma/2}.$ Then, there exists $C>0$ such that
\begin{eqnarray*}\Vert S(t)f\Vert_{L^{\infty}}&\leq& C \Big\{ \epsilon+(1-\sigma)|t|^{-1/2}\epsilon^{-1/2}\\
&& +|t|^{-1/2}\max\{N^{3/2},\epsilon^{-1/2}\}\left[ N^{\sigma}+(\sqrt{3}+\epsilon)^{\sigma}\right]-\frac{N^{\sigma+1}}{\sigma+1}\Big\}\Vert f\Vert_{L^1},
 \end{eqnarray*}
for all $f\in \mathscr{S}(\mathbb{R})$ and $t\neq 0.$
\end{lemm}

\proof
Define $h_{x,t}(\xi)=x\xi-t\varphi(\xi)$ and $F(\xi)=(1+\xi^2)^{\sigma/2}.$ Then
\begin{eqnarray*} 
K(t)f(x)&=&\int_{-\infty}^{\infty}F(\xi)e^{ih_{x,t}(\xi)}\widehat{f}(\xi)d\xi\\
&=&\int_{|\xi|\leq \epsilon}+\int_{\epsilon<|\xi|<\sqrt{3}-\epsilon}+\int_{\sqrt{3}-\epsilon\leq|\xi|\leq\sqrt{3}+\epsilon}+\int_{\sqrt{3}+\epsilon<|\xi|< N}+\int_{|\xi|\geq N}\\
&=:&I_1+I_2+I_3+I_4+I_5.
\end{eqnarray*}
From the choice of $\sigma,$ we have that  $|F(\xi)|\leq 1$ for all $\xi.$ Then using the Riemann-Lebesgue's Lemma, we obtain
\begin{equation}\label{DesI_1}
|I_1|\leq \Vert \widehat{f}\Vert_{\infty}\int_{-\epsilon}^{\epsilon}|F(\xi)|d\xi\leq 2\epsilon \Vert f\Vert_{L^1}.
\end{equation}
In a similar way we get
\begin{equation}\label{DesI_3}
|I_3|\leq 2\epsilon \Vert f\Vert_{L^1}.
\end{equation}
Next, to estimate $I_2$ note that $h_{x,t}$ is concave in $[\epsilon, \sqrt{3}-\epsilon];$ furthermore,
\[h'_{x,t}(\xi)=x-\frac{t(1-\xi^2)}{(1+\xi^2)^2}\ \ \ \ \ \  \mbox{and} \ \ \ \ \ \  h''_{x,t}(\xi)=\frac{2t\xi(3-\xi^2)}{(1+\xi^2)^3}.\]
Then, easily  we can see that 
\begin{equation}\label{Corput2A}
\vert h''_{x,t}(\xi)\vert\geq \frac{2|t|(\sqrt{3}-\epsilon)(3-(\sqrt{3}-\epsilon)^2)}{(1+(\sqrt{3}-\epsilon)^2)^3}\apprge |t| \epsilon,
\end{equation}
for all $\xi\in [\epsilon, \sqrt{3}-\epsilon].$ On the other hand, since $\sigma<0,$ we get
\begin{equation}\label{Corput3A}
|F(\sqrt{3}-\epsilon)|\leq (1+(\sqrt{3}-\epsilon)^2)^{\sigma/2}\apprle 1.
\end{equation}
Since $\sigma<0,$  for all $\xi\in [\epsilon, \sqrt{3}-\epsilon],$  we have 
\begin{equation}\label{Corput4A}
|F'(\xi)|=-\sigma|\xi|(1+\xi^2)^{\sigma/2-1}\leq -\sigma(1+\xi^2)^{\sigma/2-1/2}\leq -\sigma.
\end{equation}
Now, applying Lemma \ref{VanCorput}, from (\ref{Corput2A}), (\ref{Corput3A}) and (\ref{Corput4A})  we arrive at
\begin{eqnarray*}
\left|\int_{\epsilon}^{\sqrt{3}-\epsilon}F(\xi)e^{ih_{x,t}(\xi)}d\xi\right|&\apprle &  (1-\sigma)|t|^{-1/2}\epsilon^{-1/2}.
\end{eqnarray*}
Therefore,
\begin{equation}\label{DesI_2}
|I_2|\apprle  (2-\sigma)|t|^{-1/2}\epsilon^{-1/2}\Vert f\Vert_{L^1}.
\end{equation}
Next, in order to estimate $I_4,$ note that $h_{x,t}$ is concave in $[\sqrt{3}+\epsilon, N].$ Then,
\begin{equation}\label{Corput2}
\vert h''_{x,t}(\xi)\vert\geq \min \left\{\frac{2|t|N(3-N^2)}{(1+N^2)^3}, \frac{2|t|(\sqrt{3}+\epsilon)|3-(\sqrt{3}+\epsilon)^2|}{(1+(\sqrt{3}+\epsilon)^2)^3}\right\}\apprge |t|\min\{N^{-3},\epsilon\},
\end{equation}
for all $\xi\in [\sqrt{3}+\epsilon, N].$ On the other hand, 
\begin{equation}\label{Corput3}
|F(N)|\leq (1+N^2)^{\sigma/2}\apprle N^{\sigma},
\end{equation}
and since $\sigma$ is negative, for all $\xi\in [\sqrt{3}+\epsilon, N]$ we obtain 
\begin{equation}\label{Corput4}
|F'(\xi)|=(1+\xi^2)^{\sigma/2-1}\vert \sigma\xi\vert \leq -\sigma(1+\xi^2)^{\sigma/2-1/2}\leq -\sigma\xi^{\sigma-1}.
\end{equation}
Now, applying Lemma \ref{VanCorput}, from (\ref{Corput2}), (\ref{Corput3}) and (\ref{Corput4})  we arrive at
\begin{eqnarray}
\left|\int_{\sqrt{3}+\epsilon}^{N}F(\xi)e^{ih_{x,t}(\xi)}d\xi\right|&\apprle &|t|^{-1/2}\max\{N^{3/2},\epsilon^{-1/2}\}\left[ N^{\sigma}-\sigma\int_{\sqrt{3}+\epsilon}^{N}\xi^{\sigma-1} d\xi \right]\nonumber\\
&\apprle &|t|^{-1/2}\max\{N^{3/2},\epsilon^{-1/2}\}\left[ N^{\sigma}+(\sqrt{3}+\epsilon)^{\sigma}-N^{\sigma}\right]\nonumber\\
&\apprle &|t|^{-1/2}\max\{N^{3/2},\epsilon^{-1/2}\}\left[ N^{\sigma}+(\sqrt{3}+\epsilon)^{\sigma}\right].
\end{eqnarray}
Therefore,
\begin{equation}\label{DesI_4}
|I_4|\apprle  \Vert f\Vert_{L^1}|t|^{-1/2}\max\{N^{3/2},\epsilon^{-1/2}\}\left[ N^{\sigma}+(\sqrt{3}+\epsilon)^{\sigma}\right].
\end{equation}
To estimate $I_5,$  since $\sigma<0,$ we arrive at
\begin{eqnarray*}
|I_5|&\leq& \Vert \widehat{f}\Vert_{L^\infty}\int_{N}^{\infty}|F(\xi)|d\xi\apprle\Vert f\Vert_{L^1}\int_{N}^{\infty}(1+\xi^2)^{\frac{\sigma}{2}}d\xi \apprle \Vert f\Vert_{L^1}\int_{N}^{\infty}\xi^{\sigma}d\xi.
\end{eqnarray*}
Thus,
\begin{equation}\label{DesI_5}
|I_5|\apprle  -\frac{N^{\sigma+1}}{\sigma+1} \Vert f\Vert_{L^1}.
\end{equation}
From (\ref{DesI_1}),  (\ref{DesI_3}), (\ref{DesI_2}), (\ref{DesI_4}) and (\ref{DesI_5}), we obtain the desired result.
\endproof

In a similar way as in Lemma \ref{EstEpsEne}, we also obtain the following result.
\begin{lemm}\label{LemmHsigmap}
Let  $\sigma<-1$, $2\leq p \leq  \infty$ and $1/p+1/{p^{\prime }}=1.$ There is $C_{\sigma}>0$ such that
\[\Vert S(t)f\Vert_{H_p^{\sigma}}\leq C_{\sigma} \vert t\vert^{2\left(\frac 12-\frac{1}{p}\right)\beta_{\sigma}}\Vert f\Vert_{L^{p'}},\]
for all $f\in \mathscr{S}(\mathbb{R})$ and $t\neq 0.$ Here $\beta_{\sigma}$ is given by (\ref{Beta}).
\end{lemm}
\proof
From Lemma \ref{EstEpsEne}, taking $N^{3/2}=\epsilon^{-1/2},$ we obtain
\begin{eqnarray*}
\Vert K(t)f\Vert_{L^{\infty}}\!&\!\apprle\!&\! \left\{ N^{-3}\!+\!(1-\sigma)|t|^{-1/2}N^{3/2}\! +\!|t|^{-1/2}N^{3/2}\left[ \!N^{\sigma}+(\sqrt{3}\!+\!N^{-3})^{\sigma}\right]-\frac{N^{\sigma+1}}{\sigma+1} \right\}\Vert f\Vert_{L^1}\\
\!&\!\apprle\!&\! \left\{ N^{-3}+(1-\sigma)|t|^{-1/2}N^{3/2} +|t|^{-1/2}N^{3/2}\left[ N^{\sigma}+1\right]-\frac{N^{\sigma+1}}{\sigma+1}  \right\}\Vert f\Vert_{L^1}\\
\!&\!\apprle\!&\! \left\{ N^{-3}+(1-\sigma)|t|^{-1/2}N^{3/2} +|t|^{-1/2}N^{\sigma+3/2}+ |t|^{-1/2}N^{3/2}-\frac{N^{\sigma+1}}{\sigma+1}  \right\}\Vert f\Vert_{L^1}\\
\!&\!\apprle\!&\! \left\{ N^{-3}+(2-\sigma)|t|^{-1/2}N^{3/2} +|t|^{-1/2}N^{\sigma+3/2}-\frac{N^{\sigma+1}}{\sigma+1}  \right\}\Vert f\Vert_{L^1}.
\end{eqnarray*}
Let $t>1,$ $\theta >0$ and $N=2t^{\theta}.$ Then, $N>2$ and
\begin{eqnarray*}
\Vert K(t)f\Vert_{L^{\infty}}&\apprle& \left\{ t^{-3\theta}+(2-\sigma)t^{\frac{3\theta}{2} -\frac 12} +t^{\theta(\sigma+\frac 32)-\frac 12}-\frac{t^{(\sigma+1)\theta}}{\sigma+1} \right\}\Vert f\Vert_{L^1}.
\end{eqnarray*}
Taking $\theta=\frac{1}{1-2\sigma},$ with $\sigma<-1,$ we obtain
\begin{eqnarray*}
\Vert K(t)f\Vert_{L^{\infty}}\leq C_{\sigma} t^{\beta_{\sigma}}\Vert f\Vert_{L^1}.
\end{eqnarray*}
Note that $C_{\sigma}\to +\infty$ as $\sigma\to -1^-$ or $\sigma\to -\infty.$ If $0<t\leq 1$ and $\sigma<-1,$ one easily get
\begin{eqnarray*}
\vert K(t)f(x)\vert&=& \left\vert\int_{-\infty}^{\infty}e^{ix\xi-it\varphi(\xi)}(1+\xi^2)^{\sigma/2}\widehat{f}(\xi)d\xi\right\vert\apprle \Vert f\Vert_{L^1}\int_{-\infty}^{\infty}(1+\xi^2)^{\frac{\sigma}{2}}d\xi\\
&\apprle& t^{\beta_{\sigma}}\Vert f\Vert_{L^1}.
\end{eqnarray*}
On the other hand, it is clear that $K(t):L^2\longrightarrow L^2$ is continuous; thus an interpolation argument and recalling that $S(t)=K(t)J^\sigma,$ permit us to finishes the proof of the lemma. 
\endproof

Next lemma gives a time-decay estimate of the group $S(t)$ on modulation spaces $M^s_{p,q}.$
\begin{lemm}\label{GrupoMs} 
Let $ s \in \mathbb{R},$ $\sigma<-1,$ $2\leq p <\infty,$ $\frac{1}{p}+\frac{1}{p'}=1,$ $0<q<\infty,$ $\theta\in[0,1],$ $\beta_{\sigma}$ as in (\ref{Beta}). Then we have
\begin{equation}
\Vert S(t) f \Vert_{M^{s}_{p,q}}\apprle (1+|t|)^{2\theta(\frac12-\frac1p)\beta_{\sigma}}\Vert f\Vert_{M^{s-\sigma\theta}_{p',q}}.
\end{equation}
\end{lemm}
\proof From Lemma \ref{LemmHsigmap}, and taking into account that $S(t)$ and $\square_k$ commutate, we obtain 
\begin{equation}\label{IneSqut}
\Vert \square_k S(t)f \Vert_{H^{\sigma}_p}\apprle \vert t\vert^{2\left(\frac 12-\frac{1}{p}\right)\beta_{\sigma}}\Vert \square_k f\Vert_{L^{p'}}. 
\end{equation}
Using the Berstein's multiplier estimate (cf. Wang \cite{Wang}), we have
\begin{equation}\label{IneBerMul}
\Vert \square_k (I-\Delta)^{\delta/2}g \Vert_{L^p}\apprle (1+ \vert k\vert)^{\delta}\Vert g\Vert_{L^{p}}.
\end{equation}
Then, from (\ref{IneSqut}) and (\ref{IneBerMul}), we arrive at
\begin{equation}\label{SquSum1}
\Vert \square_k S(t) f \Vert_{L^p} \apprle (1+ \vert k\vert)^{-\sigma}\sum_{l\in\Lambda} \Vert \square_{k+l} S(t) f \Vert_{H^{\sigma}_p} \apprle (1+ \vert k\vert)^{-\sigma} \vert t\vert^{2\left(\frac 12-\frac{1}{p}\right)\beta_{\sigma}} \sum_{l\in\Lambda} \Vert \square_{k+l} f \Vert_{L^{p'}}.
\end{equation}
On the other hand, from the H\"older and Young's inequalities, we obtain
\begin{align}\label{SquSum2}
\Vert \square_k S(t)f \Vert_{L^p}& \apprle  \Vert \sigma_k e^{-it\varphi(\xi)}\varphi(\xi)\widehat{f} \Vert_{L^{p'}} \apprle \sum_{l\in\Lambda}\Vert \sigma_k e^{-it\varphi(\xi)}\varphi(\xi)\mathscr{F}\square_{k+l}f \Vert_{L^{p'}}\notag\\
&\apprle \sum_{l\in\Lambda}\Vert \mathscr{F}\square_{k+l}f \Vert_{L^p} \apprle  \sum_{l\in\Lambda}\Vert \square_{k+l}f \Vert_{L^{p'}}.
\end{align}
From (\ref{SquSum1}) and (\ref{SquSum2}) and an interpolation argument we get 
\begin{equation}\label{SquTheta}
\Vert \square_k S(t) f \Vert_{L^p} \apprle (1+ \vert k\vert)^{-\sigma\theta} \vert t\vert^{2\theta\left(\frac 12-\frac{1}{p}\right)\beta_{\sigma}} \sum_{l\in\Lambda} \Vert \square_{k+l} f \Vert_{L^{p'}},
\end{equation}
for any $\theta\in[0,1].$ Since $-\sigma\geq 0,$ from (\ref{SquSum2}) we have
\begin{align}\label{SquSumThe}
\Vert \square_k S(t) f \Vert_{L^p} \apprle  (1+ \vert k\vert)^{-\sigma\theta} \sum_{l\in\Lambda}\Vert \square_{k+l}f \Vert_{L^{p'}}. 
\end{align}
Combining (\ref{SquTheta}) and (\ref{SquSumThe}), we arrive at
\begin{equation}\label{FinIneST}
\Vert \square_k S(t) f \Vert_{L^p} \apprle (1+ \vert k\vert)^{-\sigma\theta}(1+ \vert t\vert)^{2\theta\left(\frac 12-\frac{1}{p}\right)\beta_{\sigma}} \sum_{l\in\Lambda} \Vert \square_{k+l} f \Vert_{L^{p'}}.
\end{equation}
Finally, multiplying (\ref{FinIneST}) by $(1+ \vert k\vert)^{s}$ and then taking the $l^p-$norm, we obtain the desired result. 
\endproof
 \begin{prop}\label{p20}
Let $s\in\mathbb{R}, \sigma<-1,$ $2\leq p<\infty,$ $0<q<\infty,$ $1\leq r \leq\infty,$ and $\theta\in [0,1].$   Then
\begin{equation}
\Vert S(t) u_0\Vert_{L^r(I,M^s_{p,q})}\apprle \Vert u\Vert_{M^{s-\theta\sigma}_{p',q}},
\end{equation}
where $I$ is a compact subset of $\mathbb{R}$ containing zero. 
\end{prop}
\proof
From Lemma \ref{GrupoMs} and since $t \mapsto (1+\vert t \vert)^{2\theta (\frac 12-\frac 1p)\beta_{\sigma}}$ is a continuos function on $I,$ we obtain
\begin{align*}
\Vert S(t)u_0\Vert_{L^r( I,M^s_{p,q})}& \apprle \left\Vert  (1+\vert t \vert)^{2\theta(\frac 12-\frac 1p)\beta_{\sigma}} \left\Vert u_0\right\Vert_{M^{s-\theta\sigma}_{p',q}}\right\Vert_{L_t^r}\\
& =\Vert u_0 \Vert_{M^{s-\theta\sigma}_{p',q}} \left\Vert  (1+\vert t \vert)^{2\theta(\frac 12-\frac 1p)\beta_{\sigma}} \right\Vert_{L_t^r(I)} \apprle \Vert u_0 \Vert_{M^{s-\theta\sigma}_{p',q}},
\end{align*}
which concludes the result.
 \endproof
 \begin{lemm}\label{Otro_GrupoMs} 
Let $ s \in \mathbb{R},$  $2\leq p <\infty,$  $0<q<\infty,$ $\theta\in[0,1]$. Then, we have
\begin{equation}
\Vert S(t)\varphi(D) f \Vert_{M^{s}_{p,q}}\apprle \langle t\rangle^{2(\frac12-\frac1p)}\Vert f\Vert_{M^{s}_{p,q}}.
\end{equation}
Here, the simbol $\langle t\rangle$ denotes $(1+\vert t\vert^2)^{1/2}.$
\end{lemm}
\proof From Lemma 3.1 of Wang \cite{Ming_Wang} it holds
\begin{equation}\label{e1}
\Vert S(t) \varphi(D)f \Vert_{L^{p}}\apprle \langle t\rangle^{2(\frac12-\frac1p)}\Vert \varphi(D)f\Vert_{L^{p}}.
\end{equation}
Notice that since $p\geq 2$ we can write $\partial_x(1-\partial_{xx})^{-1}(f)=\partial_xJ^{\frac{1}{p}-2}J^{-\frac{1}{p}}(f),$ where $J^s$ is the fractional differential operator defined by the symbol $(1+\vert\xi \vert^2)^{s/2}.$ Since the symbol of the operator $\partial_xJ^{\frac{1}{p}-2}$ is $i\xi\langle \xi\rangle^{\frac{1}{p}-2},$ it is an $L^p-$multiplier (cf. Theorem 2.1 in \cite{Ming_Wang}). Furthermore, taking into account that the mapping $J^{-\frac{1}{p}}:L^{p/2}\rightarrow L^p$ is bounded, it holds that $\Vert \varphi(D)f \Vert_{L^{p}}\apprle \Vert f\Vert_{L^{p}}.$ Consequently, from (\ref{e1})  we get
\begin{equation}\label{e2}
\Vert S(t) \varphi(D)f \Vert_{L^{p}}\apprle \langle t\rangle^{2(\frac12-\frac1p)}\Vert f\Vert_{L^{p}}.
\end{equation}
From (\ref{e2}) and since $\square_k$ and $S(t)\varphi(D)$ commutate, we have
\begin{equation}\label{e3z}
\Vert \square_kS(t) \varphi(D)f \Vert_{L^{p}}\apprle \langle t\rangle^{2(\frac12-\frac1p)}\Vert \square_k f\Vert_{L^{p}}.
\end{equation}
Multiplying (\ref{e3z}) by $(1 + \vert k\vert)^s$ and taking the $l^q-$norm in both sides of (\ref{e3z}), we obtain the
desired result.
\endproof

\begin{lemm}\label{LemmVarphi}
Define $\varphi_k(\xi)=\varphi(\xi-k),$ where $\varphi(\xi)=\frac{\xi}{1+\xi^2}.$ Then,
\[ \left\vert \frac{\partial^2 \varphi_k(\xi)}{\partial\xi^2}\right\vert \apprle \langle k \rangle^{-3},\]
for all $|\xi|\leq C$ and $k\in\mathbb{Z}.$
\end{lemm}
\proof
Notice that 
\[ \frac{\partial^2 \varphi_k(\xi)}{\partial\xi^2}=\frac{2(\xi+k)((\xi+k)^2-3)}{(1+(\xi+k)^2)^3}.\]
Therefore,
\[ \left| \frac{\partial^2 \varphi_k(\xi)}{\partial\xi^2}\right|\leq \frac{6(1+|\xi+k|^2)^{\frac 12}(1+|\xi+k|^2)}{(1+|\xi+k|^2)^3}=6{(1+|\xi+k|^2)^{-\frac32}}.\]
Since $|\xi|\leq C,$ we obtain the desired result.
\endproof
\begin{prop}\label{Isomor}
Let $1<p\leq \infty,$  $0<q\leq \infty,$  $s \in\mathbb{R}.$ Then,
\[
\Vert \varphi(D)g\Vert_{M_{p,q}^{s}}\apprle \Vert g\Vert_{M_{p,q}^{s-1}},
\]
for all $g\in M_{p,q}^{s-1}.$
\end{prop}
\proof
The proof is inspired in the proof of Proposition A.1 in Kato \cite{Kato}. First, we choose an auxiliary smooth function $\kappa\in\mathscr{S}$ satisfying
\begin{equation*}
\kappa(\xi)=\left\{
\begin{array}{lc}
1, & \mbox{if} \ \ \vert\xi\vert\leq 1,\vspace{0.2cm}\\
0, & \mbox{if} \ \ \vert\xi\vert \geq 2.
\end{array}
\right. 
\end{equation*}
We also define
\[\kappa_k(\xi):=\kappa\left(\frac{\xi-k}{C}\right).\]
Then $\kappa_k=1$ on the support of $\sigma_k.$ Here, the constant $C>1$ is that one taking from the support of $\sigma_k$ in Section \ref{intro}.  By the Young's inequality and the  change of variables $\xi-k\mapsto \xi$ we have
\begin{align*}
\Big\Vert \mathscr{F}^{-1}[\sigma_k \cdot \mathscr{F}(\varphi(D)g)]\Big\Vert_{L^p}&=\Big\Vert \mathscr{F}^{-1}[ \kappa_k \varphi\cdot \sigma_k  \mathscr{F}g]\Big\Vert_{L^p}\\
&\leq \Big\Vert \mathscr{F}^{-1}[ \kappa_k \varphi]\Big\Vert_{L^1} \Big\Vert \mathscr{F}^{-1}[\sigma_k  \mathscr{F}g]\Big\Vert_{L^p}\\
&= \Big\Vert \int_{\mathbb{R}} e^{ix\xi}\kappa\left(\frac{\xi-k}{C}\right) \varphi(\xi) d\xi\Big\Vert_{L^1} \Vert \mathscr{F}^{-1}[\sigma_k  \mathscr{F}g]\Vert_{L^p}\\
&= \left\Vert \int_{\mathbb{R}} e^{ix\xi}\kappa\left(\frac{\xi}{C}\right) \varphi(\xi-k) d\xi\right\Vert_{L^1} \Big\Vert \square_k g\Big\Vert_{L^p}.
\end{align*}
Now, since $|\varphi(\xi-k)|\apprle \langle k\rangle^{-1}$ for all $|\xi|\leq C$ and for all $k\in \mathbb{Z},$ we obtain
\begin{align*}
 \int_{\vert x \vert \leq 1}\left| \int_{\mathbb{R}} e^{ix\xi}\kappa\left(\frac{\xi}{C}\right) \varphi(\xi-k) d\xi\right| dx\apprle  \int_{\vert x \vert \leq 1} \langle k\rangle^{-1} dx\apprle  \langle k\rangle^{-1}.
 \end{align*}
Notice that
 \[
 \frac{1}{x}\frac{\partial e^{ix\xi}}{\partial \xi}=ie^{ix\xi}, \ \ \ \text{for} \ \ \ x\neq 0.
 \]
 Then, integrating by parts twice, from Lemma \ref{LemmVarphi} we easily see  that
 \begin{align*}
 \int_{\vert x \vert \geq  1}\left| \int_{\mathbb{R}} e^{ix\xi}\kappa\left(\frac{\xi}{C}\right) \varphi(\xi-k) d\xi\right| dx\apprle  \int_{\vert x \vert \geq 1} \frac{1}{x^2} \langle k\rangle^{-1} dx\apprle  \langle k\rangle^{-1}.
 \end{align*}
 Thus
\[
\Big\Vert \mathscr{F}^{-1}[\sigma_k \cdot \mathscr{F}(\varphi(D)g)]\Big\Vert_{L^p}\apprle \langle k\rangle^{-1} \Vert \square_k g\Vert_{L^p}.
\]
Now, with the usual modifications in the case $q=\infty$, we have 
 \[
 \Vert \varphi(D)g\Vert_{M_{p,q}^{s}}=\left( \sum_{k\in\mathbb{Z}}\langle k\rangle^{sq} \Vert \square_k \varphi(D)g\Vert^q_{L^p} \right)^{1/q}\apprle \left( \sum_{k\in\mathbb{Z}}\langle k\rangle^{(s-1)q}\Vert \square_k g\Vert^q_{L^p} \right)^{1/q} \apprle  \Vert g\Vert_{M_{p,q}^{s-1}},
 \]
 as desired.
\endproof

With the aim of making the reading easier, we present three lemmas which allow us to deal with the nonlinearity $f(u)=u^{\lambda+1}$. The proof of the first two can be found in Iwabuchi \cite{Iwabuchi} (Proposition 2.7 (ii) and Corollary 2.9 (ii)). For the proof of the third one we refer to B\'enyi and Okoudjou \cite{Okoudjou}.
 \begin{lemm} \label{Ib1}
 Let $1\leq p,p_1,p_2\leq \infty,$ $1<\sigma,\sigma_1,\sigma_2<\infty.$   If $\frac{1}{p}=\frac{1}{p_1}+\frac{1}{p_2},$ $\frac{1}{\sigma}-\frac{1}{\sigma_1}-\frac{1}{\sigma_2}+1\leq \frac{s}{n}<\frac{1}{\sigma},$ there exists $C>0,$ such that for any $u\in M^s_{p_1,\sigma_1}(\mathbb{R}^n)$ and $v\in M^s_{p_2,\sigma_2}(\mathbb{R}^n),$ it holds
 \[ \|u v\|_{M^{s}_{p,\sigma}}\leq C  \|u\|_{M^{s}_{p_1,\sigma_1}}\|v\|_{M^{s}_{p_2,\sigma_2}}.\]
 \end{lemm}

\begin{lemm}\label{Ib2}
Let $1\leq q \leq \infty, $ $p \in \mathbb{N},$ $0\leq  s <n/{\nu},$ and $1\leq \mu, \nu <\infty$ satisfying
\[ \frac{1}{\nu}-\frac{(p-1)s}{n}\leq \frac{p}{\mu}-p+1, \ \ \ \ \ 1\leq \nu\leq \mu.\]
Then, there exists $C>0$ such that for any $u\in M^{s}_{{pq},\mu}(\mathbb{R}^n),$ we have
\[ \|u^p\|_{M^{s}_{{q},\mu}}\leq C  \|u\|^p_{M^{s}_{{pq},\nu}}.\]
\end{lemm}

\begin{lemm}\label{prod}(Product estimate) Let $m$ a positive integer and $s\geq 0.$ Assume that $\frac{1}{p_1}+\cdots +\frac{1}{p_m}=\frac{1}{p_0},$ $\frac{1}{q_1}+\cdots +\frac{1}{q_m}=m-1+\frac{1}{q_0},$ with $0<p_i\leq \infty,$ $1\leq q_i\leq \infty$ for $1\leq i\leq m.$ Then it holds that
\[ \left\| \prod_{i=1}^mu_i\right\|_{M^{s}_{{p_0},q_0}}\leq C \prod_{i=1}^m \|u_i\|^p_{M^{s}_{{p_i},q_i}},\]
where $C$ is independent of $u_i$, and $u_i\in M^{s}_{{p_i},q_i},$ $i=1,...,m.$
\end{lemm}

 \begin{prop}\label{propuse}
Let $\lambda\geq 1$ an integer, $p=\lambda+2,$ $1\leq q <\infty,$ $1-\frac 1q \leq s<\frac 1q,$ $r=\frac{\lambda(\lambda+2)}{\lambda+2+\theta\lambda \beta_{\sigma}},$    and $f(u)=u^{\lambda+1}.$  Also assume that $0<\theta\leq -\frac{1}{\sigma}$ and $\sigma<-1$ such that $0<\frac{-\theta\lambda\beta_{\sigma}}{\lambda+2}<1.$ Then
\begin{equation}
\left\Vert \int_0^t S(t-\tau)\varphi(D)f(u(\tau))d\tau\right\Vert_{L^r(\mathbb{R},M^s_{p,q})}\apprle \Vert u\Vert^{\lambda+1}_{L^r(\mathbb{R},M^s_{p,q})}.
\end{equation}
\end{prop}
\proof
From Lemma \ref{GrupoMs} and Proposition \ref{Isomor} we obtain
\begin{align*}
\left\Vert \int_0^t S(t-\tau)\varphi(D)f(u(\tau))d\tau\right\Vert_{L^r(\mathbb{R},M^s_{p,q})}&\apprle \left\Vert  \int_0^t \left\Vert S(t-\tau)\varphi(D)f(u(\tau))\right\Vert_{M^s_{p,q}} d\tau\right\Vert_{L_t^r}\\
& \apprle \left\Vert  \int_0^t (1+\vert t-\tau\vert)^{2\theta(\frac 12-\frac 1p)\beta_{\sigma}} \left\Vert \varphi(D)f(u(\tau))\right\Vert_{M^{s-\sigma\theta}_{p',q}} d\tau\right\Vert_{L_t^r}\\
& \apprle \left\Vert  \int_0^t (1+\vert t-\tau\vert)^{2\theta(\frac 12-\frac 1p)\beta_{\sigma}} \left\Vert f(u(\tau))\right\Vert_{M^{s-\sigma\theta-1}_{p',q}} d\tau\right\Vert_{L_t^r}.
\end{align*}
Since $\sigma<-1,$ we have $-\frac{1}{\sigma}<1.$ Therefore, we can choose $0<\theta\leq -\frac{1}{\sigma},$ and apply the embedding $M^s_{p',q}\subset M^{s-\sigma\theta-1}_{p',q}$ and Lemma \ref{Ib2} in last inequality to arrive at
\begin{align}
\left\Vert \int_0^t S(t-\tau)\varphi(D)f(u(\tau))d\tau\right\Vert_{L^r(\mathbb{R},M^s_{p,q})}&\apprle  \left\Vert  \int_0^t (1+\vert t-\tau\vert)^{2\theta(\frac 12-\frac 1p)\beta_{\sigma}} \left\Vert f(u(\tau))\right\Vert_{M^{s}_{p',q}} d\tau\right\Vert_{L_t^r}\nonumber\\
&\apprle  \left\Vert  \int_0^t (1+\vert t-\tau\vert)^{2\theta(\frac 12-\frac 1p)\beta_{\sigma}} \Vert u(\tau)\Vert^{\lambda+1}_{M^{s}_{p,q}} d\tau\right\Vert_{L_t^r}.\label{z12}
\end{align}
From the choice of $\theta, \sigma$ and $r,$ we have that
\[ 0<1+2\theta\left( \frac 12-\frac 1p\right)\beta_{\sigma}<1 \ \ \ \  \text{and}\ \ \ \  \frac 1 r=\frac{\lambda+1}{r}-\left(1+2\theta\left( \frac 12-\frac 1p\right)\beta_{\sigma}\right).\]
Therefore, we can apply the Hardy-Littlewood-Sobolev's inequality in (\ref{z12}) in order to obtain the desired result.
 \endproof
 
 \begin{prop}\label{alter}
Let $\lambda\geq 1$ an integer, $p=\lambda+2,$ $s\geq 0,$ $r=\frac{\lambda(\lambda+2)}{\lambda+2+\theta\lambda \beta_{\sigma}},$  and $f(u)=u^{\lambda+1}.$  Also assume that $0<\theta\leq -\frac{1}{\sigma}$ and $\sigma<-1$ such that $0<\frac{-\theta\lambda\beta_{\sigma}}{\lambda+2}<1.$ Then
\begin{equation}
\left\Vert \int_0^t S(t-\tau)\varphi(D)f(u(\tau))d\tau\right\Vert_{L^r(\mathbb{R},M^s_{p,1})}\apprle \Vert u\Vert^{\lambda+1}_{L^r(\mathbb{R},M^s_{p,1})}.
\end{equation}
\end{prop}
\proof
The proof of Proposition \ref{alter} is analogous to the proof of Proposition \ref{propuse} by using Lemma \ref{prod} in place of Lemma \ref{Ib2}
\endproof

\section{General Strichartz estimates in $M^{s}_{p,q}$}\label{Section3}
The aim of this section is to derive some Strichartz estimates in modulation spaces $M^s_{p,q}$ for a general dispersive semigroup 
\begin{equation}\label{Ufmenosqpf}
U(t)=\mathscr{F}^{-1}e^{itP(\xi)}\mathscr{F},
\end{equation}
where $P:\mathbb{R}^n\to\mathbb{R}$ is a real-valued function. In Section 4, we will use these estimates in the particular case of $U(t)=S(t)$ in order to analyze the existence of global solutions for the gBBM equation. We assume that $U(t)$ is a semigroup which satisfies the next estimate
\begin{equation}\label{Umpq}
 \Vert U(t)f\Vert_{M^{\alpha}_{p,q}}\apprle (1+|t|)^{-\mu}  \Vert f\Vert_{M^{\alpha+\delta}_{p',q}},
 \end{equation}
where $2\leq p <\infty,$ $1\leq q<\infty,$ $\alpha \in \mathbb{R},$ $\mu=\mu(p),$ with  $0<\mu<1$ and $\alpha, \delta, \mu$ are independent of $t\in \mathbb{R}.$ Taking into account  Lemma \ref{GrupoMs}, it holds that an example of a group $U(t)$ verifying  (\ref{Ufmenosqpf}) and (\ref{Umpq}) is given by the BBM group $S(t).$
Next, we derive some time-space estimates of $U(t)$ satisfying (\ref{Ufmenosqpf}) and (\ref{Umpq}). Here is worthwhile to remark that the results of this section can be applied to analyze existence results for other dispersive models in the framework of modulation spaces.  These estimates complement those proved by Wang and Hudzik \cite{BaoxHudz}. We need to establish some additional notations: Given the Banach space $X,$ for instance $X=L^\gamma(\mathbb{R};L^p(\mathbb{R}^n)),$ $1\leq p,\gamma\leq \infty,$ we consider the function spaces $l_{\square}^{s,q}(X),$ $s\in\mathbb{R},1\leq q<\infty,$ 
introduced in \cite{BaoxHudz}, which are defined as follows:
\begin{align*}
l_{\square}^{s,q}(X)=\left\{ f\in \mathscr{S}^{\prime}(\mathbb{R}^{n+1}): \Vert u\Vert_{l_{\square}^{s,q}(X)}:= \left( \sum_{k\in\mathbb{Z}^n}(1+\vert k\vert)^{sq}\Vert \square_kf\Vert^q_X \right)^{1/q}<\infty\right\}. 
\end{align*}
In the definition of $l_{\square}^{s,q}(X)$, if $s = 0,$ then we write $l_{\square}^{q}(X).$ The following duality results is known (cf. Wang and Hudzik \cite{BaoxHudz}).
\begin{theo}(Dual space)\cite{BaoxHudz}\label{DualSpace}
 Let $s\in \mathbb{R}$ and $1\leq p,\gamma<\infty.$ We have
\[(c_{\square}^s(L^{\gamma}(\mathbb{R},L^p)))^*=l_{\square}^{-s,1}(L^{\gamma'}(\mathbb{R},L^{p'})).\]
\end{theo}
\begin{prop}\label{Strichartz1}
Let $U(t)$ satisfying (\ref{Ufmenosqpf}) and (\ref{Umpq}). Then, for all $s\in\mathbb{R}$ and $\gamma=\frac{2}{\mu}$ we have
\begin{equation}\label{Strilm2q}
\Vert U(t)f\Vert_{l^{s,q}_{\square}(L^{\gamma}(\mathbb{R},L^p))}\apprle \Vert f\Vert_{M^{s+\frac{\delta}{2}}_{2,q}}.
\end{equation}
In addition, if $\gamma\geq q,$ we have
\begin{equation}\label{StriLM2q}
\Vert U(t)f\Vert_{L^{\gamma}(\mathbb{R}, M^{s}_{p,q})}\apprle \Vert f\Vert_{M^{s+\frac{\delta}{2}}_{2,q}}.
\end{equation}
\end{prop}
\proof The proof is based on a duality argument. Without loss of generality take $s=-\delta/2.$ First, we consider the case $1<q<\infty.$ We show that
\begin{equation}\label{IntM2qlpri}
\int_{\mathbb{R}}(U(t)f,\phi(t))dt\apprle \Vert f\Vert_{M_{2,q}}\Vert \phi \Vert_{l_{\square}^{\delta/2,q'}(L^{\gamma'}(\mathbb{R},L^{p'}))}
\end{equation}
holds for all $f\in \mathscr{S}(\mathbb{R}^n)$ and $\phi\in C_0^{\infty}(\mathbb{R}; \mathscr{S}(\mathbb{R}^n)).$ Since  $\mathscr{S}(\mathbb{R}^n)$ and $C_0^{\infty}(\mathbb{R};\mathscr{S}(\mathbb{R}^n))$ are dense in $M_{2,q}$ and  $l_{\square}^{\delta/2,q'}(L^{\gamma'}(\mathbb{R},L^{p'})),$ respectively (cf. \cite{BaoxHudz}), we see that (\ref{IntM2qlpri}) implies (\ref{Strilm2q}). By duality we obtain that
\begin{equation}\label{IntM2qM2qpri}
\int_{\mathbb{R}}(U(t)f,\phi(t))dt\apprle \Vert f\Vert_{M_{2,q}}\left\Vert \int_{\mathbb{R}}U(-t)\phi(t)dt \right\Vert_{M_{2,q'}}.
\end{equation}
Now, for $k\in \mathbb{Z}^n$ we get
\begin{equation}\label{LgammaPri}
\left\Vert \square_k \int_{\mathbb{R}}U(-t)\phi(t)dt \right\Vert^2_{L^2} \apprle \Vert \square_k \phi\Vert_{L^{\gamma'}(\mathbb{R};L^{p'})}\left\Vert \square_k\int_{\mathbb{R}}U(t-\tau)\phi(\tau)d\tau \right\Vert_{L^{\gamma}(\mathbb{R};L^p)}.
\end{equation}
Since the sequence $\{\square_k\}_{k\in\mathbb{Z}}$ is almost orthogonal, using (\ref{Umpq}) with $\alpha=-\delta$, the definition of the norm of $M_{p',q}$ and the Bernstein's multiplier estimate (cf. Wang and Huang \cite{Wang}) we arrive at
\begin{align}\label{Uktkpri}
\Vert \square_k U(t)f \Vert_{L^p}&\apprle(1+\vert k\vert)^{\delta} \Vert U(t)f \Vert_{M^{-\delta}_{p,q}}\apprle (1+\vert t\vert)^{-\mu(p)}(1+\vert k\vert)^{\delta} \sum_{l\in\Lambda}\Vert \square_k \square_{k+l}f \Vert_{M_{p',q}}\nonumber\\
&\apprle (1+\vert t\vert)^{-\mu(p)}(1+\vert k\vert)^{\delta} \Vert \square_k f \Vert_{L^{p'}},
\end{align}
where $\Lambda=\{l\in \mathbb{Z}^n: B(0,\sqrt{n})\cup B(l,\sqrt{n})\neq \emptyset\}.$ Since $0<\mu<1$ and $\gamma=\frac{2}{\mu}$ we can use the Hardy-Littlewood-Sobolev's inequality to obtain
\begin{equation}\label{CuaIntPhi}
\left\Vert \square_k \int_{\mathbb{R}}U(t-\tau)\phi(\tau)d\tau \right\Vert_{L^{\gamma}(\mathbb{R};L^p)} \apprle (1+\vert k\vert)^{\delta} \Vert \square_k \phi\Vert_{L^{\gamma'}(\mathbb{R};L^{p'})}.
\end{equation}
So, in view of (\ref{LgammaPri}) and (\ref{CuaIntPhi}), we have
\[
\left\Vert \square_k \int_{\mathbb{R}}U(-t)\phi(t)dt \right\Vert_{L^2} \apprle (1+\vert k\vert)^{\delta/2} \Vert \square_k \phi\Vert_{L^{\gamma'}(\mathbb{R};L^{p'})}.
\]
Taking the $l^{q'}-$norm in both sides of the last inequality, we get
\begin{equation}\label{IntUlpeq}
\left\Vert \int_{\mathbb{R}}U(-t)\phi(t)dt \right\Vert_{M_{2,q'}} \apprle  \Vert \phi\Vert_{l_{\square}^{\delta/2,q'}(L^{\gamma'}(\mathbb{R};L^{p'})}.
\end{equation}
From (\ref{IntM2qlpri}) and (\ref{IntUlpeq}), we arrived at (\ref{Strilm2q}), as desired.\\

If $\gamma\geq q,$ using the Minkowski's inequality, we obtain the left-hand of (\ref{Strilm2q}) is controlled by the left-hand of (\ref{IntM2qlpri}).\\

Next, we consider the case $q=1.$ From Theorem \ref{DualSpace}, it enough to show that 
\[
\int_{\mathbb{R}}(U(t)f,\phi(t))dt\apprle \Vert f\Vert_{M_{2,q}}\Vert \phi \Vert_{c_{\square}^{\delta/2}(L^{\gamma'}(\mathbb{R},L^{p'}))}
\]
holds for all $f\in \mathscr{S}(\mathbb{R}^n)$ and $\phi\in C_0^{\infty}(\mathbb{R}; \mathscr{S}(\mathbb{R}^n)).$ Repeating the above procedure, we can obtain the desire result. Finally, from the Minkowski's inequality we obtain (\ref{StriLM2q}) from  (\ref{Strilm2q}).
\endproof

Next, we estimate the nonlinear part. We denote by
\[(\mathcal{N} f)(t,\cdot)=\int_0^tU(t-\tau)f(\tau, \cdot)d\tau.\]
\begin{prop}\label{Strichartz2}
 Let $U(t)$ satisfying (\ref{Ufmenosqpf}) and (\ref{Umpq}). Then, for all $s\in\mathbb{R}$ and $\gamma=\frac{2}{\mu}$ we have
\begin{equation}\label{Strilqamed}
\Vert \mathcal{N} f\Vert_{l^{s,q}_{\square}(L^{\infty}(\mathbb{R},L^2))}\apprle \Vert f\Vert_{l^{s+\frac{\delta}{2},q}_{\square}(L^{\gamma'}(\mathbb{R},L^{p'}))}.
\end{equation}
In addition, if $\gamma'\leq  q,$ we have
\begin{equation}\label{StriLinfgam}
\Vert \mathcal{N} f\Vert_{L^{\infty}(\mathbb{R}, M^s_{2,q})}\apprle \Vert f\Vert_{L^{\gamma'}(\mathbb{R},M^{s+\frac{\delta}{2}}_{p',q})}.
\end{equation}
\end{prop}
\proof 
From the definition of the norm of the space $l^{s,q}_{\square}(L^{\infty}(\mathbb{R},L^2))$ and using the same ideas as in Proposition \ref{Strichartz1}, the crucial inequality 
\[
\Vert \square_k \mathcal{N}f \Vert^2_{L^2}\apprle \langle k\rangle^{\delta} \Vert \square_k f \Vert^2_{L^{\gamma'}(\mathbb{R};L^{p'})},
\]
implies (\ref{Strilqamed}). From the Minkowski's inequality we obtain (\ref{StriLinfgam}) from  (\ref{Strilqamed}). 

\begin{prop}
 Let $U(t)$ satisfying (\ref{Ufmenosqpf}) and (\ref{Umpq}). Then, for  all $s\in\mathbb{R}$ and $\gamma=\frac{2}{\mu},$  we have
\begin{equation}\label{StrilalqS}
\Vert \mathcal{N} f\Vert_{l^{s,q}_{\square}(L^{\gamma}(\mathbb{R},L^p))}\apprle \Vert f\Vert_{l^{s+\frac{\delta}{2},q}_{\square}(L^1(\mathbb{R},L^2))}.
\end{equation}
In addition, if $\gamma\geq q,$ we have
\begin{equation}\label{StriLgL1}
\Vert \mathcal{N} f\Vert_{L^{\gamma}(\mathbb{R}, M^{s}_{p,q})}\apprle  \Vert f\Vert_{L^1(\mathbb{R},M^{s+\frac{\delta}{2}}_{2,q})}.
\end{equation}
\end{prop}
\proof 
Let $f,\phi\in C_0^{\infty}(\mathbb{R},\mathscr{S}(\mathbb{R})).$ Following Proposition \ref{Strichartz1}, we have
\begin{align*}
\left| \int_{\mathbb{R}}\left(\int_0^t U(t-\tau)f(\tau)d\tau,\phi(t)\right) dt\right|&\apprle \Vert f\Vert_{L^1(\mathbb{R},M_{2,q})} \left\Vert \int_{\cdot}^{\infty}U(\cdot-t)\phi(t)dt\right\Vert_{L^{\infty}(\mathbb{R}; M_{2,q'})}\\
& \apprle \Vert f\Vert_{L^1(\mathbb{R},M_{2,q})} \left\Vert \phi \right\Vert_{l_{\square}^{\frac{\delta}{2},q}(L^{\gamma'}(\mathbb{R};L^{p?}))}.
\end{align*}
Since $\phi\in C_0^{\infty}(\mathbb{R},\mathscr{S}(\mathbb{R}))$ is dense in $l_{\square}^{\frac{\delta}{2},q}(L^{\gamma'}(\mathbb{R};L^{p'}))$ and in $c_{\square}^{\delta/2}(L^{\gamma'}(\mathbb{R},L^{p'})),$ by duality we get (\ref{StrilalqS}). Finally, from the Minkowski's inequality we obtain (\ref{StriLgL1}) from  (\ref{StrilalqS}).
\endproof
\begin{prop}\label{Prop1}
 Let $U(t)$ satisfying (\ref{Ufmenosqpf}) and (\ref{Umpq}). Then, for all $s\in\mathbb{R}$  and  $\gamma=\frac{2}{\mu},$  we have
\begin{equation}\label{StriLaLma}
\Vert \mathcal{N} f\Vert_{l^{s,q}_{\square}(L^{\gamma}(\mathbb{R},L^p))}\apprle \Vert f\Vert_{l^{s+\frac{\delta}{2},q}_{\square}(L^{\gamma'}(\mathbb{R},L^{p'}))}.
\end{equation}
In addition, if $q\in [\gamma',\gamma],$ we have
\begin{equation}\label{StriLgLgp}
\Vert \mathcal{N} f\Vert_{L^{\gamma}(\mathbb{R}, M^{s}_{p,q})}\apprle \Vert f\Vert_{L^{\gamma'}(\mathbb{R},M^{s+\delta}_{p',q})}.
\end{equation}
\end{prop}
\proof 
The proof of \ref{StriLaLma} is also based on a duality argument, so, we only sketch the proof of (\ref{StriLgLgp}).  From (\ref{Umpq}) we obtain
\[\Vert \mathcal{N} f\Vert_{M^{s}_{p,q}}\apprle  \int_0^t \langle t-\tau\rangle^{-\mu} \left\Vert f(\tau)\right\Vert_{M^{s+\delta}_{p',q}} d\tau.\]
 Taking the $L^{\gamma}-$norm and using the Hardy-littlewood-Sobolev inequality we obtain the desired result. 
 \endproof
\section{Existence results}
\subsection{Proof of Theorems \ref{GlobalLamb0} and \ref{GlobalLamb1}}
We first focus on the proof of Theorem \ref{GlobalLamb1}, which is based on a fixed point argument by applying the time-decay and Strichartz estimates, as well as the nonlinear estimates obtained in Section 2 and 3. For that, let us consider the closed ball $B_{2\epsilon}=\{u:\Vert u\Vert_{L^{r}(\mathbb{R}; M^s_{p,q})}\leq 2\epsilon\},$ with $\epsilon>0,$ and define the map $\Gamma $ on the metric space $B_{2\epsilon}:$
\[
(\Gamma u)(x,t)=S(t)u_{0}(x)-\frac{i}{\lambda+1}\int_{0}^{t}S(t-\tau )\varphi(D)[u^{\lambda+1}(x,\tau)] d\tau.
\]
We want to choose a $\epsilon>0$ such that $\Gamma: B_{2\epsilon}\rightarrow B_{2\epsilon}$ is a contraction. From Proposition \ref{propuse} and the smallness assumption on $\Vert S(t)u_0\Vert_{L^r(\mathbb{R};M^s_{p,q})}$ we have that, if $u\in B_{2\epsilon},$ then
\begin{align}
\Vert \Gamma u\Vert_{L^{r}(\mathbb{R};M^s_{p,q})}&\leq \Vert S(t)u_{0}(x)\Vert_{L^{r}(\mathbb{R};M^s_{p,q})}+C\Vert u\Vert^{\lambda+1}_{L^{r}(\mathbb{R}; M^s_{p,q})}\notag\\
& \leq \epsilon+C(2\epsilon)^{\lambda+1}=\epsilon+C2^{\lambda+1}\epsilon^\lambda\epsilon.\label{e18}
\end{align}
Taking $\epsilon>0$ such that $C2^{\lambda+1}\epsilon^\lambda<1,$ we get that $\Gamma:B_{2\epsilon}\rightarrow B_{2\epsilon}.$ From Lemma \ref{GrupoMs} and Proposition \ref{Isomor} (see the proof of Proposition \ref{propuse}) and taking into account Lemmas \ref{Ib1} and \ref{Ib2}, we get
\begin{align*}
&\Vert \Gamma u-\Gamma v\Vert_{L^{r}(\mathbb{R};M^s_{p,q})}\apprle  \left\Vert  \int_0^t (1+\vert t-\tau\vert)^{2\theta(\frac 12-\frac 1p)\beta_{\sigma}} \left\Vert u^{\lambda+1}-v^{\lambda+1}\right\Vert_{M^{s-\sigma\theta}_{p',q}} d\tau\right\Vert_{L_t^r}\\
&\ \ \ \ \apprle  \left\Vert  \int_0^t (1+\vert t-\tau\vert)^{2\theta(\frac 12-\frac 1p)\beta_{\sigma}} \left\Vert (u-v)\left( \sum_{k=1}^{\lambda+1}u^{\lambda+1-k} v^{k-1}\right)\right\Vert_{M^{s-\sigma\theta}_{p',q}} d\tau\right\Vert_{L_t^r}\\
&\ \ \ \ \apprle  \left\Vert  \int_0^t (1+\vert t-\tau\vert)^{2\theta(\frac 12-\frac 1p)\beta_{\sigma}} \left\Vert u-v\right\Vert_{M^s_{p,q}}\sum_{k=1}^{\lambda+1}\Vert u^{\lambda+1-k}\Vert_{M^s_{\frac{p}{\lambda+1-k},q}}\Vert v^{k-1}\Vert_{M^s_{\frac{p}{k-1},q}}d\tau\right\Vert_{L_t^r}\\
&\ \ \ \ \apprle  \left\Vert  \int_0^t (1+\vert t-\tau\vert)^{2\theta(\frac 12-\frac 1p)\beta_{\sigma}} \Vert u-v\Vert_{M^s_{p,q}}\sum_{k=1}^{\lambda+1}\Vert u\Vert^{\lambda+1-k}_{M^s_{p,q}}\Vert v\Vert^{k-1}_{M^s_{p,q}}d\tau\right\Vert_{L_t^r}.
\end{align*}
Considering the assumption $0<\frac{-\lambda\theta\beta_{\sigma}}{\lambda+2}<1$ and $r=\frac{\lambda(\lambda+2)}{\lambda+2+\lambda\theta\beta_{\sigma}},$ we can  apply the H\"older and Hardy-littlewood-Sobolev's inequalities in last inequality in order to obtain
\begin{align}
\Vert \Gamma u-\Gamma v\Vert_{L^{r}(\mathbb{R};M^s_{p,q})}&\leq \Vert u-v\Vert_{L^{r}(\mathbb{R};M^s_{p,q})}(\Vert u\Vert^\lambda_{L^{r}(\mathbb{R};M^s_{p,q})}+\Vert v\Vert^\lambda_{L^{r}(\mathbb{R};M^s_{p,q})})\nonumber\\
&\leq 2^{\lambda+1}\epsilon^\lambda C\Vert u-v\Vert_{L^{r}(\mathbb{R};M^s_{p,q}).}\label{e19}
\end{align}
From (\ref{e18}) and (\ref{e19}) it holds that $\Gamma:B_{2\epsilon}\rightarrow B_{2\epsilon}$ is a contraction, which implies the existence of a unique fixed point, as desired. The proof of Theorem \ref{GlobalLamb0} is analogous to the proof of Theorem \ref{GlobalLamb1} by using Proposition \ref{alter} in place of Proposition \ref{propuse}.
\endproof
\subsection{Proof of Theorem \ref{TheoX0}}
The proof is based on a fixed point argument by applying the time-decay and Strichartz estiamtes, as well as the nonlinear estimates obtained in Section 2 and Section 3. For that, let us consider the closed ball $B_{2\epsilon}=\{u:\Vert u\Vert_{L^{\infty}(\mathbb{R}; M^s_{2,q})}+\Vert u\Vert_{L^{\gamma}(\mathbb{R}; M^s_{p,q})}\leq 2\epsilon\},$ with  $\epsilon>0,$
 and define the map $\Gamma $ on the metric space $B_{2\epsilon}:$ 
\[
(\Gamma u)(x,t)=S(t)u_{0}(x)-\frac{i}{\lambda+1}\int_{0}^{t}S(t-\tau )\varphi(D)[u^{\lambda+1}(x,\tau)] d\tau.
\]
From Proposition \ref{Strichartz1} (Inequality (\ref{StriLM2q})) and Proposition \ref{Prop1} (Inequality (\ref{StriLgLgp})) we have
\begin{align*}
\Vert \Gamma u\Vert_{L^{\gamma}(\mathbb{R};M^s_{p,q})}&\apprle \Vert u_0\Vert_{M^{s-\frac{\sigma\theta}{2}}_{2,q}}+\Vert \varphi(D)(u^{\lambda+1})\Vert_{L^{\gamma'}(\mathbb{R};M^{s-\sigma\theta}_{p',q})}.
\end{align*}
In view of Proposition \ref{Isomor}, the choice of $\theta,$ Lemma \ref{Ib2} and the values of $p$ and $\gamma,$ we arrive at
\begin{align}\label{EstLgamma}
\Vert \Gamma u\Vert_{L^{\gamma}(\mathbb{R};M^s_{p,q})}&\apprle \Vert u_0\Vert_{M^{s-\frac{\sigma\theta}{2}}_{2,q}}+\Vert u^{\lambda+1}\Vert_{L^{\gamma'}(\mathbb{R};M^{s-\sigma\theta-1}_{p',q})}\notag\\
&\apprle \Vert u_0\Vert_{M^{s-\frac{\sigma\theta}{2}}_{2,q}}+\Vert u^{\lambda+1}\Vert_{L^{\gamma'}(\mathbb{R};M^{s}_{p',q})}\notag\\
&\apprle \Vert u_0\Vert_{M^{s-\frac{\sigma\theta}{2}}_{2,q}}+\Vert u\Vert^{\lambda+1}_{L^{\gamma'(\lambda+1)}(\mathbb{R};M^{s}_{p'(\lambda+1),q})}\notag\\
&\apprle \Vert u_0\Vert_{M^{s-\frac{\sigma\theta}{2}}_{2,q}}+\Vert u\Vert^{\lambda+1}_{L^{\gamma}(\mathbb{R};M^{s}_{p,q})}.
\end{align}
On the other hand, since  $\Vert \square_k S(t)u_0\Vert_2\leq \Vert \square_k u_0\Vert_2 $ and Proposition \ref{Strichartz2} (Inequality (\ref{StriLinfgam})), we have
\begin{align*}
\Vert \Gamma u\Vert_{L^{\infty}(\mathbb{R};M^s_{2,q})}&\apprle \Vert u_0\Vert_{M^{s}_{2,q}}+\Vert \varphi(D)(u^{\lambda+1})\Vert_{L^{\gamma'}(\mathbb{R};M^{s-\frac{\sigma\theta}{2}}_{p',q})}.
\end{align*}
In view Proposition \ref{Isomor}, the choice of $\theta,$ Lemma \ref{Ib2} and the values of $p$ and $\gamma,$ we obtain
\begin{align}\label{EstLinfty}
\Vert \Gamma u\Vert_{L^{\infty}(\mathbb{R};M^s_{2,q})}&\apprle \Vert u_0\Vert_{M^{s}_{2,q}}+\Vert u^{\lambda+1}\Vert_{L^{\gamma'}(\mathbb{R};M^{s-\frac{\sigma\theta}{2}-1}_{p',q})}\notag\\
&\apprle \Vert u_0\Vert_{M^{s-\frac{\sigma\theta}{2}}_{2,q}}+\Vert u^{\lambda+1}\Vert_{L^{\gamma'}(\mathbb{R};M^{s}_{p',q})}\notag\\
&\apprle \Vert u_0\Vert_{M^{s-\frac{\sigma\theta}{2}}_{2,q}}+\Vert u\Vert^{\lambda+1}_{L^{\gamma}(\mathbb{R};M^{s}_{p,q})}.
\end{align}
From (\ref{EstLgamma}) and (\ref{EstLinfty}), if $\Vert u_0\Vert_{M^{s-\frac{\sigma\theta}{2}}_{2,q}}\leq \epsilon$ and $u\in B_{2\epsilon},$ we conclude that $\Gamma:B_{2\epsilon}\rightarrow B_{2\epsilon}.$ Also, in the same spirit of the proof of Theorem \ref{GlobalLamb1} together with (\ref{EstLgamma}) and (\ref{EstLinfty}) we get that $\Gamma:B_{2\epsilon}\rightarrow B_{2\epsilon}$ is a contraction on $B_{2\epsilon}$, which implies the existence of a unique fixed point, as desired. Now we prove that the fixed point $u\in C(\mathbb{R};M^s_{2,q}),$ that is, that $\Vert u(t_n)-u(t)\Vert_{M^s_{2,q}}\rightarrow 0,$ as $t_n\rightarrow t.$ For that, notice that
\begin{eqnarray}
\Vert u(t_n)-u(t)\Vert_{M^s_{2,q}} &\apprle & \Vert S(t_n)u_0-S(t)u_0\Vert_{M^s_{2,q}}+\frac{1}{\lambda+1}\Big\Vert \int_0^{t_n}S(t_n-\tau)\varphi(D)[u^{\lambda+1}(x,\tau)] d\tau\nonumber\\
& &-\int_0^{t}S(t-\tau)\varphi(D)[u^{\lambda+1}(x,\tau)] d\tau\Big\Vert_{M^s_{2,q}}\nonumber\\
 &:= & I_1+\frac{1}{\lambda+1}I_2.\label{lm1}
\end{eqnarray}
Observe that for $v\in \mathscr{S}^\Omega=\{f:f\in \mathscr{S}\ \mbox{and}\ \mbox{supp}\widehat{f}\subset \Omega\}$ (see Section \ref{intro}) it holds
\begin{eqnarray*}
\left\Vert \square_k(S(t_n)v-S(t)v)\right\Vert_2 &\apprle & \left \Vert \sigma_k\left (e^{-it_n\frac{\xi}{1+\xi^2}}-e^{-it\frac{\xi}{1+\xi^2}}\right)\widehat{v}(\xi)\right\Vert_2\nonumber\\
&\apprle & \left \Vert \left (e^{-it_n\frac{\xi}{1+\xi^2}}-e^{-it\frac{\xi}{1+\xi^2}}\right)\widehat{v}(\xi)\right\Vert_2.
\end{eqnarray*}
If $v\in \mathscr{S},$ then from the Lebesgue's dominated convergence Theorem we get $$\left \Vert \left (e^{-it_n\frac{\xi}{1+\xi^2}}-e^{-it\frac{\xi}{1+\xi^2}}\right)\widehat{u_0}(\xi)\right\Vert_2\rightarrow 0,\ \mbox{as}\ t_n\rightarrow t.$$ 
Since $v\in \mathscr{S}^\Omega,$ then $\square_k(S(t_n)v-S(t)v)\neq 0$ only for a finite number of $k.$ Therefore we can conclude that
$$I_1=\Vert S(t_n)u_0-S(t)u_0\Vert_{M^s_{2,q}}\rightarrow 0, \ \mbox{as}\ t_n\rightarrow t.$$
Now we deal with the term $I_2$ in (\ref{lm1}). For that, notice that
\begin{eqnarray*}
I_2&= &\Big\Vert (S(t_n)-S(t))\int_0^t S(-\tau)\varphi(D)[u^{\lambda+1}(x,\tau)] d\tau+\int_t^{t_n}S(t_n-\tau)\varphi(D)[u^{\lambda+1}(x,\tau)] d\tau\Big\Vert_{M^s_{2,q}}\\
&:= & I_2^1+I_2^2.
\end{eqnarray*}
Thus, arguing  as in (\ref{EstLinfty}) we get
\begin{eqnarray*}
\Big\Vert \int_0^{t}S(-\tau)\varphi(D)[u^{\lambda+1}(x,\tau)] d\tau\Big\Vert_{M^s_{2,q}} &\apprle & \left\Vert u\right\Vert^{\lambda+1}_{L^\gamma(\mathbb{R};M^s_{p,q})}<\infty,
\end{eqnarray*}
which allow us to conclude that $I_2^1\rightarrow 0,\ \mbox{as}\ t_n\rightarrow t.$ In the same way, 
\begin{eqnarray*}
\Big\Vert \int_t^{t_n}S(t_n-\tau)\varphi(D)[u^{\lambda+1}(x,\tau)] d\tau\Big\Vert_{M^s_{2,q}} &\apprle & \left\Vert u\right\Vert^{\lambda+1}_{L^\gamma(\mathbb{R};M^s_{p,q})}<\infty,
\end{eqnarray*}
which implies that $I_2^2\rightarrow 0,\ \mbox{as}\ t_n\rightarrow t.$ Thus, $I_2\rightarrow 0,\ \mbox{as}\ t_n\rightarrow t,$ and then $u\in C(\mathbb{R};M^s_{2,q}).$ 
\endproof


\subsection{Proof of Theorems \ref{TheoX2a} and \ref{TheoX2}}
We first present the proof of Theorem \ref{TheoX2}. This proof is also based on a fixed point argument and the time-decay and product estimates established in Section 2.  Let us consider the closed ball 
\[B_{2\epsilon}=\left\{u:\sup_{-\infty<t<\infty}(1+\vert t\vert)^{\rho}\Vert u(t)\Vert_{M^s_{p,q}}\leq 2\epsilon\right\},\ \epsilon>0,\]
where $\rho=-2\theta(\frac{1}{2}-\frac{1}{p})\beta_\sigma>0,$ and define the map $\Gamma $ on the metric space $B_{2\epsilon}:$
 \[
(\Gamma u)(x,t)=S(t)u_{0}(x)-\frac{i}{\lambda+1}\int_{0}^{t}S(t-\tau )\varphi(D)[u^{\lambda+1}(x,\tau)] d\tau.
\]
From Lemma \ref{GrupoMs}, Proposition \ref{Isomor}, the embedding $M^{s}_{p',q}\subset M^{s-{\sigma\theta}-1}_{p',q},$ and Lemma \ref{Ib2}, we have 
\begin{align}\label{EstLinf}
\Vert \Gamma u\Vert_{M^s_{p,q}}&\apprle (1+\vert t\vert)^{-\rho}\Vert u_0\Vert_{M^{s-{\sigma\theta}}_{p',q}}+\int_0^t (1+\vert t-\tau\vert)^{-\rho} \Vert \varphi(D)u^{\lambda+1}\Vert_{M^{s-{\sigma\theta}}_{p',q}}\notag\\
&\apprle (1+\vert t\vert)^{-\rho}\Vert u_0\Vert_{M^{s-{\sigma\theta}}_{p',q}}+\int_0^t (1+\vert t-\tau\vert)^{-\rho}\Vert u^{\lambda+1}\Vert_{M^{s-{\sigma\theta}-1}_{p',q}}\notag\\
&\apprle (1+\vert t\vert)^{-\rho}\Vert u_0\Vert_{M^{s-{\sigma\theta}}_{p',q}}+\int_0^t (1+\vert t-\tau\vert)^{-\rho}\Vert u^{\lambda+1}\Vert_{M^{s}_{p',q}}\notag\\
&\apprle (1+\vert t\vert)^{-\rho}\Vert u_0\Vert_{M^{s-{\sigma\theta}}_{p',q}}+\int_0^t (1+\vert t-\tau\vert)^{-\rho}\Vert u\Vert^{\lambda+1}_{M^{s}_{p,q}}\notag\\
&\apprle (1+\vert t\vert)^{-\rho}\Vert u_0\Vert_{M^{s-{\sigma\theta}}_{p',q}}+\Big(\sup_{t>0}(1+\vert t\vert)^{\rho}\Vert u(t)\Vert_{M^s_{p,q}}\Big)^{\lambda+1}\notag\\
&\times \int_0^t (1+\vert t-\tau\vert)^{-\rho}(1+\vert \tau\vert)^{-\rho(\lambda+1)}d\tau.
\end{align}
Since $\lambda\geq 3,$ then $\rho(\lambda+1)>1.$ Thus, it holds that
\begin{align}
\int_0^{t/2}(1+t-\tau)^{-\rho}(1+ \tau)^{-\rho(\lambda+1)}d\tau \apprle (1+t)^{-\rho}\int_0^{t/2}(1+ \tau)^{-\rho(\lambda+1)}d\tau\apprle (1+ t)^{-\rho}.\label{mod2qno1}
\end{align}
Also, it is straightforward to get
\begin{align}
\int_{t/2}^t(1+ t-\tau)^{-\rho}(1+ \tau)^{-\rho(\lambda+1)}d\tau&= \int_{0}^{t/2}(1+ \tau)^{-\rho}(1+ t-\tau)^{-\rho(\lambda+1)}d\tau \nonumber\\
&\ \apprle(1+t)^{-\rho(\lambda+1)}\int_0^{t/2}(1+ \tau)^{-\rho}d\tau  \apprle (1+ t)^{-\rho}.\label{mod3qno1}
\end{align}
Therefore, from (\ref{EstLinf})-(\ref{mod3qno1}) we obtain
\begin{align}\label{EstLinf2}
\sup_{t>0}(1+\vert t\vert)^{\rho}\Vert \Gamma u\Vert_{M^s_{p,q}}&\apprle \Vert u_0\Vert_{M^{s-{\sigma\theta}}_{p',q}}+\Big(\sup_{t>0}(1+\vert t\vert)^{\rho}\Vert u(t)\Vert_{M^s_{p,q}}\Big)^{\lambda+1}.
\end{align}
From (\ref{EstLinf2}), if $\Vert u_0\Vert_{M^{s-{\sigma\theta}}_{p',q}}\leq \epsilon$ and $u\in B_{2\epsilon},$ we conclude that $\Gamma:B_{2\epsilon}\rightarrow B_{2\epsilon}.$ Also, following the proof of (\ref{EstLinf}) we get that $\Gamma:B_{2\epsilon}\rightarrow B_{2\epsilon}$ is a contraction on $B_{2\epsilon}$, which implies the existence of a unique fixed point, as desired. The proof of Theorem \ref{TheoX2a} follows analogously to the proof of Theorem \ref{TheoX2} by using Lemma \ref{prod} in place of Lemma \ref{Ib2}.
\endproof
\subsection{Proof of Theorems \ref{TheoX10} and \ref{TheoX10b}} 
We first prove Theorem \ref{TheoX10b}. The proof of Theorem \ref{TheoX10b} is also based on a fixed point argument. Let us consider the closed ball 
\[
B_{2\epsilon}=\left\{u\in C([0,T];M^s_{p,q}):\ \sup_{0<t<T}\Vert u(t)\Vert_{M^s_{p,q}}\leq 2\epsilon\right\},\ \epsilon>0,
\]
 and define the map $\Gamma $ on the metric space $B_{2\epsilon}$ \[
(\Gamma u)(x,t)=S(t)u_{0}(x)-\frac{i}{\lambda+1}\int_{0}^{t}S(t-\tau )\varphi(D)[u^{\lambda+1}(x,\tau)] d\tau.
\]
From Lemma \ref{GrupoMs} with $\theta =0,$ and Lemma \ref{Otro_GrupoMs}, the embedding $M^s_{\frac{p}{\lambda+1},q}\subset M^s_{p,q}$ and Lemma \ref{Ib2} we get
\begin{align}
\Vert \Gamma u\Vert_{M^s_{p,q}}& \apprle  (1+\vert t\vert)^{2\theta(\frac12-\frac1p)\beta_{\sigma}} \Vert u_0\Vert_{M^{s-\sigma\theta}_{p',q}}+\int_0^t\Vert S(t-\tau)\varphi(D)u^{\lambda+1}(\tau)\Vert_{M^s_{p,q}}d\tau\nonumber\\
& \apprle  (1+\vert t\vert)^{2\theta(\frac12-\frac1p)\beta_{\sigma}} \Vert u_0\Vert_{M^{s-\sigma\theta}_{p',q}}+\int_0^t  \langle t-\tau \rangle^{2(\frac 12-\frac 1p)}\Vert u^{\lambda+1}(\tau)\Vert_{M^s_{p,q}}d\tau\nonumber\\
& \apprle  (1+\vert t\vert)^{2\theta(\frac12-\frac1p)\beta_{\sigma}} \Vert u_0\Vert_{M^{s-\sigma\theta}_{p',q}}+\int_0^t  \langle t-\tau \rangle^{2(\frac 12-\frac 1p)}\Vert u^{\lambda+1}(\tau)\Vert_{M^s_{\frac{p}{\lambda+1},q}}d\tau\nonumber\\
& \apprle  (1+\vert t\vert)^{2\theta(\frac12-\frac1p)\beta_{\sigma}} \Vert u_0\Vert_{M^{s-\sigma\theta}_{p',q}}+\int_0^t  \langle t-\tau \rangle^{2(\frac 12-\frac 1p)}\Vert u(\tau)\Vert^{\lambda+1}_{M^s_{p,q}}d\tau. \label{e3}
\end{align}
Taking $\theta=0$ in (\ref{e3}) we get
\begin{align}
\Vert \Gamma u\Vert_{M^s_{p,q}}& \leq 
C\Vert u_0\Vert_{M^{s}_{p',q}}+C\langle T \rangle^{2(\frac 12-\frac 1p)}T\Vert u\Vert^{\lambda+1}_{M^s_{p,q}}. \label{e4}
\end{align}
Let $\epsilon=C\Vert u_0\Vert_{M^s_{p',q}}$ and consider $T>0$ such that $C\langle T \rangle^{2(\frac 12-\frac 1p)}T2^{\lambda+1}\epsilon^\lambda<1.$ Then, if we assume that $\sup\limits_{0<t<T}\Vert u(t)\Vert_{M^s_{p,q}}\leq 2\epsilon,$ from (\ref{e4}) we get
\begin{align}\label{e21}
\Vert \Gamma u\Vert_{M^s_{p,q}} \leq 
\epsilon+C\langle T \rangle^{2(\frac 12-\frac 1p)}T(2\epsilon)^{\lambda+1}\leq 2\epsilon. 
\end{align}
Therefore, $\Gamma:{B}_{2\epsilon}\longrightarrow {B}_{2\epsilon}.$ In an analogous way we get
\begin{align}
\Vert \Gamma u-\Gamma v\Vert_{M^s_{p,q}}&\leq C\int_0^t  \langle t-\tau \rangle^{2(\frac 12-\frac 1p)}\Vert u^{\lambda+1}(\tau)-v^{\lambda+1}(\tau)\Vert_{M^s_{p,q}}d\tau\nonumber\\
&\leq C\sup_{0<t<T}\Vert u-v\Vert_{M^s_{p,q}}\langle T \rangle^{2(\frac 12-\frac 1p)}T2^{\lambda+1}\epsilon^\lambda,
\end{align}
which implies that $\Gamma$ is a contraction on $B_{2\epsilon}$, and thus, the integral equation (\ref{IntEqu}) has a unique local solution $u\in B_{2\epsilon}.$ The proof of the time-continuity follows in the same way to the one of Theorem \ref{TheoX0}, and therefore we omit it. The proof of Theorem \ref{TheoX10} follows analogously to the proof of Theorem \ref{TheoX10b}  by using Lemma \ref{prod} in place of Lemma \ref{Ib2}.

\endproof



\begin{thebibliography}{99}

\bibitem{AngBanSci1} J. Angulo, C.  Banquet and  M. Scialom, \textit{Stability for the modified and fourth-order Benjamin-Bona-Mahony equations,} Discrete  Continuous Dynamical Systems - A, 30 (2011), 851-871. 

\bibitem{AngBanSci2} J. Angulo, C.  Banquet and M. Scialom, \textit{The regularized Benjamin-Ono and BBM equations: Well-posedness and nonlinear stability,} J. Differential Equations, 250 (2011), 4011-4036. 

\bibitem{Avrin} J. Avrin and J.  Goldstein,  \textit{Global existence for the Benjamin-Bona-Mahony equation in arbitrary dimensions,} Nonlinear Analysis, 9 (1985), 861-865.

\bibitem{BenBonMah} T. Benjamin, J. Bona and J. Mahony, \textit{Model equations for long waves in nonlinear dispersive systems,} Phil.Trans. R. Soc., 272 (1972), 47-78.

\bibitem{Okoudjou} \'A.B\'enyi and K. Okoudjou, \textit{Local well-posedness of nonlinear dispersive equations on modulation spaces,} Bulletin of the London Mathematical Society, 41, (2009), 549-558.
 
\bibitem{BonCas} J. Bona and R. Cascaval, \textit{Nonlinear dispersive waves on trees,} Can. Appl. Math. Q., 16 (2008), 1-18.

\bibitem{BonDai} J. Bona and M Dai, \textit{Norm-inflation results for the BBM equation,} J. Math. Anal. Appli., 446 (2017), 879 - 885.

\bibitem{Bona} J.  Bona and  N. Tzvetkov, \textit{Sharp well-posedness results for the BBM equation,} Discrete Contin. Dyn. Syst, 23 (2009), 1241-1252.

\bibitem{BonPriSco} J. Bona, W. Pritchard and  L. Scott, \textit{A comparison of solutions of two model equations for long waves,} Lectures in Applied Mathematics, 20 (1983), 235-267.

\bibitem{Carvajal} X. Carvajal and M. Panthee,  \textit{On ill-posedness for the generalized BBM equation,} Discrete and Contin. Dyn. Syst., 34  (2014), 4565-4576.

\bibitem{Chai} L. Chaichenets, D. Hundertmark, P. Kunstmann and  N. Pattakos,  \textit{On the existence of global solutions of the one-dimensional cubic NLS for initial data in the modulation space $M_{p,q}(\mathbb{R}),$} J. Differential Equations, 263 (2017) 4429-4441.

\bibitem{eldika}  K. El Dika, \textit{Asymptotic stability of solitary waves for the Benjamin-Bona-Mahony equation}, Discrete Contin. Dyn. Syst., 13 (2005),  583-622.


\bibitem{LinaresPonce} F. Linares and G. Ponce, \textit{Introduction to nonlinear dispersive equations}, Springer, New York, 2000.

\bibitem{Feichtinger} H.  Feichtinger, \textit{Modulation Spaces on Locally Compact Abelian Groups,} Technical Report, University of Vienna, 1983.

\bibitem{Goldstein} J. Goldstein and B.  Wichnoski, \textit{On the Benjamin-Bona-Mahony equation in higher dimensions,} Nonlinear Analysis, 4 (1980), 665-675.

\bibitem{haragus1}  M. H\u{a}r\u{a}gu\c{s}, \textit{Stability of periodic waves for the generalized BBM equation}, Rev. Roumaine Math. Pures Appl., 53 (2008),  445-463.

\bibitem{Huang}Q. Huang, D. Fan and J. Chen, \textit{Critical exponent for evolution equations in modulation spaces,} J. Math. Anal. Appl., 443 (2016), 230-242.

\bibitem{Iwabuchi}T.  Iwabuchi, \textit{Navier-Stokes equations and nonlinear heat equations in modulation spaces with negative derivative indices,} J. Differential Equations, 248 (2010), 1972-2002.

\bibitem{Kato} T. Kato, \textit{The inclusion relations between $\alpha$-modulation spaces and Lp-Sobolev spaces or local Hardy spaces,} Journal of Functional Analysis, 272 (2017), 1340-1405.

\bibitem{Manna} R. Manna, \textit{Modulation spaces and non-linear Hartree type equations,} Nonlinear Analysis, 162 (2017) 76-90.

\bibitem{MillerWeinstein} J.  Miller and M.  Weinstein, \textit{Asymptotic stability of solitary waves for the regularized long-wave equation}, Comm. Pure Appl. Math.,  495 (1996),  399-441.

\bibitem{Panthee}  M. Panthee, \textit{On the ill-posedness result for the BBM equation,} Discrete Contin. Dyn. Syst.,
30 (2011), 253-259.

\bibitem{Peregrine1} D.  Peregrine, \textit{Calculations of the development of an undular bore},  J. Fluid Mech.,25 (1966),  321--330.

\bibitem{Peregrine2} D.  Peregrine, \textit{Long waves on a beach}, J. Fluid Mech.,  27 (1967),  815--827.
 
\bibitem{Ruz} M. Ruzhansky, M. Sugimoto and B. Wang, \textit{Modulation spaces and nonlinear evolution equations, in: Evolution Equations of Hyperbolic and Schr\"odinger Type,} in: Progress in Mathematics, 301 (2012), 267-283.

\bibitem{Roum} D. Roum\'egoux, \textit{A symplectic non-squeezing theorem for BBM equation,} Dyn. Partial Differ. Equ., 7 (4) (2010), 289-305.

\bibitem{SougaStrauss}  P.  Souganidis and W.  Strauss, \textit{Instability of a class of dispersive solitary waves}, Proc. Roy. Soc. Edinburgh Sect. A, 114 (1990),  195-212.

\bibitem{Stein} E. Stein, \textit{Harmonic Analysis: Real-Variable Method, Orthogonality, and Oscillatory Integrals}, Princeton, NJ, 1993.

\bibitem{Wadati1}   M. Wadati, \textit{Wave propagation in nonlinear lattice I},  J. Phys. Soc. Japan, 38 (1975),  673-680.

\bibitem{Wadati2}  M. Wadati, \textit{Wave propagation in nonlinear lattice II}, J. Phys. Soc. Japan, 38 (1975),  681-686.

\bibitem{BaoxHudz}B. Wang and H. Hudzik, \textit{The global Cauchy problem for the NLS and NLKG with small rough data,} J. Differential Equations, 232 (2007), 36-73.

\bibitem{Wang2} M. Wang, \textit{Long time behavior of a damped generalized BBM equation in low regularity spaces,}
Mathematical Methods in the Applied Sciences, 38 (2015), 4852-4866.

\bibitem{Ming_Wang} M. Wang, \textit{Sharp global well-posedness of the BBM equation in $L^p$ type Sobolev spaces,} Discrete and Continuous Dynamical Systems-A 36 (2016), 5763-5788.

\bibitem{Wang} B. Wang and C. Huang, \textit{Frequency-uniform decomposition method for the generalized BO, KdV and NLS equations,} J.
Differential Equations, 239  (2007), 213-250.


\bibitem{weinstain1} M. Weinstein, \textit{Existence and dynamic stability of solitary wave solutions of equations arising in long wave propagation,} Comm. PDE, 12 (1987),  1133-1173.

\bibitem{Whitham} G. Whitham, \textit{Linear and Nonlinear Waves, Wiley,} New York, 1974.

\bibitem{zeng} L. Zeng, \textit{Existence and stability of solitary-wave solutions of equations of Benjamin-Bona-Mahony type}, J. Differential Equations, 188 (2003),  1-32.

\bibitem{Zhao} G. Zhao, J. Chen and W. Guo, \textit{Klein-Gordon Equations on Modulation Spaces,} Abstract and Applied Analysis, (2014), Article ID 947642, 1-15.
\end{thebibliography}
\end{document}